\documentclass[hidelinks=false,onefignum,onetabnum]{siamart250211}

\usepackage{lipsum}
\usepackage{amsfonts}
\usepackage{graphicx}
\usepackage{epstopdf}
\ifpdf
  \DeclareGraphicsExtensions{.eps,.pdf,.png,.jpg}
\else
  \DeclareGraphicsExtensions{.eps}
\fi

\newsiamremark{remark}{Remark}
\newsiamremark{hypothesis}{Hypothesis}
\crefname{hypothesis}{Hypothesis}{Hypotheses}
\newsiamthm{claim}{Claim}

\headers{HaTT: Hadamard avoiding TT recompression}{Z. Sun, J. Huang, C. Xiao and C. Yang}

\title{HaTT: Hadamard avoiding TT recompression and its applications
}

\author{Zhonghao Sun\thanks{SKLMS, AMSS, Chinese Academy of Sciences, Beijing, China, 100190  
  (\email{sunzhonghao@amss.ac.cn}).}
\and Jizu Huang\thanks{Corresponding author. SKLMS, AMSS, Chinese Academy of Sciences, Beijing, China, 100190  
  (\email{huangjz@lsec.cc.ac.cn}).}
\and Chuanfu Xiao\thanks{School of Mathematical Sciences, Peking University, and PKU-Changsha Institute for Computing and Digital Economy (\email{chuanfuxiao@pku.edu.cn})}
\and Chao Yang\thanks{School of Mathematical Sciences, Peking University, and PKU-Changsha Institute for Computing and Digital Economy (\email{chao\_yang@pku.edu.cn})}}

\usepackage{amsopn}
\DeclareMathOperator{\diag}{diag}

\ifpdf
\hypersetup{
  pdftitle={
  HaTT: Hadamard avoiding TT recompression and its applications in nonlinear problems},
  pdfauthor={Z. Sun, J. Huang, C. Xiao and C. Yang}
}
\fi

\usepackage{multirow}
\usepackage[T1]{fontenc}
\usepackage{tikz-network}

\usepackage{enumitem}
\usepackage{amssymb}
\usepackage{bm}
\usepackage[noEnd=false]{algpseudocodex}
\algrenewcommand\textproc{\textup}
\makeatletter
\algdef{SE}[FUNCTION]{Function}{EndFunction}[3]{%
	\algpx@startCodeCommand\algpx@startIndent\algorithmicfunction\ \ifstrempty{#1}{}{#1 = }\textproc{#2}\ifstrempty{#3}{}{(#3)}%
}{%
	\algpx@startEndBlockCommand\algpx@endIndent\algorithmicend\ \algorithmicfunction%
}
\algnewcommand{\LineComment}[1]{\Statex\textcolor{\algpx@commentColor}{\(\triangleright\) \itshape #1}}
\makeatother
\newcommand{\tensor}[1]{\bm{\mathcal{#1}}}
\newcommand{\mat}[1]{\bm{#1}}
\newcommand{\vb}[1]{\bm{#1}}
\newcommand{\ttcore}[2]{\tensor{T}_{\tensor{#1}, #2}}
\newcommand{\R}{\mathbb{R}}
\newcommand{\norm}[1]{\left\| #1 \right\|}
\newcommand{\fnorm}[1]{\left\| #1 \right\|_{\mathrm{F}}}

\newcommand{\alg}[1]{\relax\ifmmode\text{#1}\else #1\fi}
\newcommand{\tran}{^{\mathsf{T}}}
\newcommand{\pkp}{\boxtimes}
\DeclareMathOperator{\rank}{rank}
\DeclareMathOperator{\Range}{range}
\NewDocumentCommand\HT{m}{\mathcal{H}\left(#1\right)}    
\NewDocumentCommand\VT{m}{\mathcal{V}\left(#1\right)}    
\NewDocumentCommand\unfold{m m}{\mat{#1}_{\langle #2 \rangle}}    
\newcommand{\dd}{\,\mathrm{d}}
\usepackage[section]{placeins}
\usepackage{dsfont}
\newsiamthm{property}{Property}
\definecolor{green}{RGB} {4, 215, 106}
\begin{document}

\maketitle

\begin{abstract}
The Hadamard product of tensor train (TT) tensors is a fundamental nonlinear operation in scientific computing and data analysis. However, due to its tendency to significantly increase TT ranks, the Hadamard product poses a major computational challenge in TT tensor-based algorithms. To address this, it is crucial to develop recompression algorithms that mitigate the effects of this rank increase. Existing recompression algorithms require an explicit representation of the Hadamard product, resulting in high computational and storage costs. 
In this work, we propose a Hadamard avoiding TT recompression (HaTT) algorithm, which reduces both computational complexity and storage requirements. By leveraging the structure of the Hadamard product in TT tensors and exploiting its Hadamard product-free property, the HaTT algorithm achieves significantly lower complexity compared to existing TT recompression methods. This is confirmed through both complexity analysis and numerical experiments. Furthermore, the HaTT algorithm is applied to solve the Allen--Cahn equation, achieving substantial speedup over existing TT recompression algorithms without sacrificing accuracy.

\end{abstract}



\section{Introduction}

Tensor train (TT) decomposition, also known as matrix product states (MPS), is a powerful tool for dimension reduction with widespread applications in scientific computing \cite{cichocki2016tensor,dolgov2015polynomial,holtz2012alternating,khoromskij2018tensor,richter2021solving}, machine learning \cite{han2018unsupervised,sidiropoulos2017tensor,sozykin2022ttopt,su2020convolutional,yang2017tensor}, and quantum computing \cite{ganahl2017continuous,iaconis2024quantum,wang2017simulations,zhou2020limits}. In TT decomposition, a tensor of order $d$ is represented as the product of $d$ core tensors, each of order no more than three. This representation ensures that the memory cost grows only linearly with $d$. Moreover, most algebraic operations on TT tensors can be efficiently performed by operating directly on the TT cores, significantly reducing computational complexity and alleviating the curse of dimensionality.  Consequently, TT representation is recognized as one of the most efficient tensor network formats for solving large-scale and high-dimensional problems \cite{cichocki2016tensor,cichocki2017tensor,lee2018fundamental,oseledets2011tensor}.

The Hadamard product is a fundamental operation in tensor-based algorithms for scientific computing and data analysis \cite{bonizzoni2016tensor,doostan2009least,espig2020iterative,kressner2017recompression,risthaus2022solving}. When applied to TT tensors, the Hadamard product can be computed using the partial Kronecker product (PKP) of TT cores, which may cause a quadratic increase in TT ranks. This rank growth significantly escalates the complexity of both the operation itself and the subsequent recompression process on the resulting tensor. For TT tensors with relatively large ranks, the Hadamard product and the following recompression step can dominate the computational cost in scientific computing and data analysis. Thus, developing efficient recompression techniques to mitigate the effects of rank increase caused by the Hadamard product is of paramount importance.


\textbf{Related work.} 
The widely used recompression algorithm for TT tensors is the \alg{TT-Rounding} algorithm \cite{oseledets2011tensor}, which consists of two main steps: orthogonalization and compression (typically using the SVD). For a TT tensor with a maximal rank $r$, the computational cost of the TT-Rounding algorithm scales as $r^3$, rendering it inefficient for TT tensors with large ranks, such as those resulting from the Hadamard product of TT tensors. To improve the performance of \alg{TT-Rounding} algorithm, several parallel approaches based on randomized sketching have been proposed \cite{al2022parallel,daas2022parallel, shi2023parallel,  xie2023error}. Additionally, to circumvent the costly orthogonalization step in \alg{TT-Rounding}, three randomized algorithms have been introduced by extending randomized low-rank matrix approximation techniques \cite{al2023randomized}. Among these, the Randomize-then-Orthogonalize (RandOrth) algorithm achieves the best speedup compared to the deterministic \alg{TT-Rounding} algorithm and is likely the state-of-the-art method for TT recompression.


\textbf{Our contributions.} 
In this work, we propose a Hadamard avoiding TT recompression (HaTT) algorithm for recompressing the Hadamard product of two TT tensors, i.e. $\tensor{A}:=\tensor{Y}\odot\tensor{Z}$. This method avoids the explicit calculation of the Hadamard product $\tensor{A}$, significantly reducing both the computational and storage costs associated with the Hadamard product and the subsequent recompression procedure compared to existing TT recompression algorithms. The HaTT algorithm proceeds as follows. First, we obtain sketching matrices by partially contracting $\tensor{A}$ with a random tensor $\tensor{R}$. By exploiting the PKP structure in TT cores of $\tensor{A}$, these sketching matrices are computed directly from TT cores of $\tensor{Y}$ and $\tensor{Z}$, bypassing the explicit computation of the TT cores of $\tensor{A}$. Next, an orthogonalization sweep is performed to construct a low rank TT tensor approximation for $\tensor{A}$. Throughout this process, the PKP structure is applied, enabling the orthogonalization sweep to be directly implemented using the TT cores of $\tensor{Y}$ and $\tensor{Z}$. In summary, the HaTT algorithm recompresses the Hadamard product $\tensor{A}$ without explicitly calculating or storing it. In terms of computational complexity, the HaTT algorithm achieves a cost reduction of approximately $r$ times compared to the RandOrth algorithm, where $r$ is the maximal rank of $\tensor{Y}$ or $\tensor{Z}$. 
The accuracy and efficiency of HaTT are validated through several numerical experiments with specific datasets. Finally, we apply the HaTT algorithm to real-world problems, such as power iteration and solving the Allen--Cahn equation. These simulations demonstrate that the HaTT algorithm achieves significant speedup compared to existing TT recompression algorithms, without sacrificing accuracy.

The remainder of this paper is organized as follows. In \cref{sec:preliminaries}, we introduce some preliminary tensor notations and review existing TT recompression algorithms for the Hadamard product. The HaTT algorithm is presented in \cref{sec:hatt}, followed by an analysis of its computational complexity in  \cref{sec:complexity-analysis}. In \cref{sec:experiments}, several benchmark simulations are conducted to demonstrate the efficiency and accuracy of the HaTT algorithm. Finally, \cref{sec:conclusion} concludes the paper.


\section{Preliminaries}\label{sec:preliminaries}

\subsection{Tensor notations and operations}

In this paper, we use the bold lowercase letter to represent a vector (e.g., $\vb{a}$), the bold uppercase letter to represent a matrix (e.g., $\mat{A}$), and the bold script letter to represent a tensor (e.g., $\tensor{A}$), respectively. For a $d$th-order tensor $ \tensor{A} \in \R^{n_{1} \times n_{2} \times \dots \times n_{d}} $,
its $ (i_{1}, i_{2}, \dots, i_{d}) $th element is denoted as $ \tensor{A}(i_{1}, i_{2}, \dots, i_{d}) $ or $ a_{i_{1}, i_{2}, \dots, i_{d}} $, and the Frobenius norm of $\tensor{A}$ is defined as 
\begin{align*}
  \fnorm{\tensor{A}} = \sqrt{\sum_{i_{1},i_2,\ldots,i_d = 1}^{n_{1},n_2,\ldots,n_d} \left| \tensor{A}(i_{1}, i_{2}, \dots, i_{d}) \right|^{2}}.
\end{align*}

We then introduce several folding and unfolding operators for vector, matrix, and tensor. For the sake of convenience, we define the multi-index as $   \overline{i_{1}i_{2}\cdots i_{d}} = i_{d} + (i_{d-1} - 1)n_{d} + \cdots + (i_{1} - 1)n_{2}\cdots n_{d}$.
The {vectorization} operator vectorizes a matrix $ \mat{A}=(a_{i,j}) \in \R^{m \times n} $ into a column vector $ \alg{vec}(\mat{A}) \in \R^{mn} $ with 
\begin{align*}
  \alg{vec}(\mat{A}) = (a_{1,1}, a_{2, 1}, \dots, a_{m, 1}, a_{1, 2}, \dots, a_{m, 2}, a_{1,3}, \dots, a_{m, n})\tran.
\end{align*}
Correspondingly, the {matricization} operator folds a vector $ \vb{x} \in \R^{mn} $ into a matrix $ \mat{X} = \vb{x}|_{n}^{m}\in \R^{m \times n} $ with $ \mat{X}(:, 1)= (x_{1}, x_{2}, \dots, x_{m})\tran $, $ \mat{X}(:, 2) = (x_{m+1}, x_{m + 2}, \dots, x_{2m})\tran $, and so on. Vectorization and matricization can be generalized to tensors. For example, the mode-$(1, 2, \dots, k)$ matricization of $\tensor{A}$ reshapes it to a matrix $ \unfold{A}{k}\in \R^{(n_{1} \cdots n_k)\times (n_{k+1}\cdots n_{d})}$ that satisfies $\unfold{A}{k}(\overline{i_{1}i_{2}\ldots i_{k}}, \overline{i_{k+1}i_{k+2}\ldots i_{d}}) = \tensor{A}(i_{1}, i_{2}, \dots, i_{d}).$
In particular, we denote $\unfold{A}{1}$ and $\unfold{A}{d-1}$ as 
$ \HT{\tensor{A}} $ and $ \VT{\tensor{A}} $, which are called horizontal and vertical matricization.

The Kronecker and Hadamard products of tensors are denoted as $\otimes$ and $\odot$, respectively. Let $ \tensor{Y}, \tensor{Z} \in \R^{n_{1} \times n_{2} \times \dots \times n_{d}} $ and $ \tensor{Q} \in \R^{m_{1} \times m_{2} \times \dots \times m_{d}} $ be $d$th-order tensors, $\tensor{Y}\otimes\tensor{Q}$ and 
$\tensor{Y}\odot\tensor{Z}$ are specifically defined as
\begin{align*}
  \tensor{A} &= \tensor{Y} \otimes \tensor{Q}, & \tensor{A}(\overline{i_{1}j_{1}}, \overline{i_{2}j_{2}}, \dots&, \overline{i_{d}j_{d}})  =\tensor{Y}(i_{1}, i_{2}, \dots, i_{d})\cdot\tensor{Q}(j_{1}, j_{2}, \dots, j_{d}),\\
  \tensor{A} & = \tensor{Y} \odot \tensor{Z}, & \tensor{A}(i_{1}, i_{2}, \dots, i_{d}) & = \tensor{Y}(i_{1}, i_{2}, \dots, i_{d}) \cdot \tensor{Z}(i_{1}, i_{2}, \dots, i_{d}). 
\end{align*}
For the Kronecker product of matrices, the following important property will be used in this paper:
\begin{align}\label{eq:kronecker_product_times_vec}
  (\mat{A} \otimes \mat{B})\vb{v} = \alg{vec}\left(\mat{B}\mat{V}\mat{A}\tran\right),
\end{align}
where $ \mat{A} \in \R^{m \times n} $, $ \mat{B} \in \R^{r \times s} $, $ \vb{v} \in \R^{ns} $, and $ \mat{V} = \vb{v}|_{n}^{s} $. An incomplete Kronecker product operation, also called the PKP, is the Kronecker product of 
tensors along specified modes \cite{lee2014fundamental}. For
two third-order tensors $\tensor{Y}\in\R^{r_1\times n\times r_2}$ and $\tensor{Z}\in\R^{s_1\times n\times s_2}$, the PKP of them along mode-1 and -3 is denoted as $\tensor{A}=\left(\tensor{Y}\pkp^{1,3}\tensor{Z}\right)\in\R^{r_1s_1\times n\times r_2s_2}$, i.e., $ \tensor{A}(:, i, :) = \tensor{Y}(:, i, :) \otimes \tensor{Z}(:, i, :) $, or elementwisely,
$\tensor{A}(\overline{\alpha_1\beta_1},i,\overline{\alpha_2\beta_2})=\tensor{Y}(\alpha_1,i,\alpha_2)\cdot\tensor{Z}(\beta_1,i,\beta_2) $.


For two tensors $\tensor{Y}\in\R^{n_1\times n_2\times\cdots\times n_d}$ and $\tensor{Z}\in\R^{m_1\times m_2\times \cdots\times m_l}$ with $n_{d-k+s}=m_s$ for all $s=1,\ldots,k$, we denote the contraction of $\tensor{Y}$ and $\tensor{Z}$ along indices $\{d-k+1,\ldots,d\}$ and $\{1,2,\ldots,k\}$ as $\tensor{A}=\left(\tensor{Y}\times^{1,\ldots,k}\tensor{Z}\right)\in\R^{n_1\times\cdots\times n_{d-k}\times m_{k+1}\times\cdots\times m_l}$, elementwisely,
    \begin{multline}
        \tensor{A}(i_1,\ldots,i_{d-k},j_{k+1},\ldots,j_{l})=\\
        \sum\limits_{\alpha_1,\ldots,\alpha_{k}=1}^{m_1,\ldots,m_{k}}\tensor{Y}(i_1,\ldots,i_{d-k},\alpha_1,\ldots,\alpha_{k})\cdot\tensor{Z}(\alpha_1,\ldots,\alpha_{k},j_{k+1},\ldots,j_l).
    \end{multline}

\subsection{Tensor train tensor and its Hadamard product} \label{subsec:tensor_train}

The TT decomposition represents $\tensor{A}\in\R^{n_1\times n_2\times\cdots\times n_d}$ as the product of $d$ tensors of order at most three, i.e.,
\begin{align}
\label{eq:tt-decom}
  \tensor{A} = \ttcore{A}{1} \times^{1} \ttcore{A}{2} \times^{1} \dots \times^{1} \ttcore{A}{d},
\end{align}
where $ \{\ttcore{A}{k} \in \R^{r_{k-1} \times n_{k} \times r_{k}} :k=1,2,\dots,d\}$ are called TT cores, and $ r_{0}=r_d=1, r_{1}, \dots,r_{d-1}$ are the TT ranks of $ \tensor{A} $. Elementwisely, \cref{eq:tt-decom} can be represented as
\begin{align}
  \label{eq:tt_decom_element}
  \tensor{A}(i_{1}, i_{2}, \dots, i_{d}) = \ttcore{A}{1}(i_{1})\ttcore{A}{2}(i_{2})\cdots\ttcore{A}{d}(i_{d}).
\end{align}
Here and in the following of this paper, we denote $\ttcore{A}{k}(i_k)$ as a simplification of $\ttcore{A}{k}(:, i_k, :)$.
For convenience, we call a tensor represented by the TT decomposition \cref{eq:tt-decom,eq:tt_decom_element} as a {TT tensor}, which can be represented by a {tensor network} diagram (see \cref{fig:tensor_network_diagrams}). 
\begin{figure}[!htpb]
  \centering
  \begin{tikzpicture}
  \Vertex[x=0,y=0,label={$ \ttcore{A}{1} $},position={above}]{X1}
  \Vertex[x=1.5,y=0,label={$ \ttcore{A}{2} $},position={above}]{X2}
  \Vertex[x=3,y=0,label={$ \ttcore{A}{3} $},position={above}]{X3}
  \Vertex[x=6,y=0,label={$ \ttcore{A}{d - 1} $},position={above}]{X4}
  \Vertex[x=7.5,y=0,label={$ \ttcore{A}{d} $},position={above}]{Xd}
  \node (Y1) at (0, -1.5) {};
  \node (Y2) at (1.5, -1.5) {};
  \node (Y3) at (3, -1.5) {};
  \node (Y4) at (6, -1.5) {};
  \node (Yd) at (7.5, -1.5) {};
  \Edge[label={$ r_{1} $}, position={below}](X1)(X2)
  \Edge[label={$ r_{2} $}, position={below}](X2)(X3)
  \Edge[label={$ r_{3} $}, position={below}](X3)(4, 0)
  \Edge(5, 0)(X4)
  \Edge[style={dashed}](4, 0)(5, 0)
  \Edge(X4)(Xd)
  \node[below, color={black!75}] at (6.75, 0) {\scriptsize $ r_{d - 1} $};
  \node[below, color={black!75}] at (5.3, 0) {\scriptsize $ r_{d - 2} $};
  \Edge[label={$ n_{1} $}, position={left}](X1)(Y1)
  \Edge[label={$ n_{2} $}, position={left}](X2)(Y2)
  \Edge[label={$ n_{3} $}, position={left}](X3)(Y3)
  \Edge[label={$ n_{d - 1} $}, position={left}](X4)(Y4)
  \Edge[label={$ n_{d} $}, position={left}](Xd)(Yd)
\end{tikzpicture}
  \caption{Tensor network diagram for a TT tensor $\tensor{A}$. }\label{fig:tensor_network_diagrams}
\end{figure}
Particularly, if the mode-$1$ (-$3$) matricization of the $k$th TT core $\ttcore{A}{k}$ satisfies row (column) orthogonality, we call $\ttcore{A}{k}$ is left (right) orthogonal.
In addition, we define the partial contracted product $ \ttcore{A}{k:l} $ for $1\leq k<l\leq d$ as
\begin{align*}
  \ttcore{A}{k:l} & = \ttcore{A}{k} \times^{1} \ttcore{A}{k + 1} \times^{1} \dots \times^{1} \ttcore{A}{l} \in \R^{r_{k-1}\times n_{k} \times \dots \times n_{l} \times r_{l}}.
\end{align*}
With TT representation, the memory cost of $\tensor{A}$ can be reduced from $  
 \prod_{k=1}^dn_k$ to $\sum_{k=1}^dn_kr_{k-1}r_k$, which grows linearly with $d$. More importantly, the basic operations of TT tensors can be converted into the corresponding operations on TT cores.


A tensor operator $ \tensor{F}: \R^{n_{1} \times \dots \times n_{d}} \to \R^{m_{1} \times \dots \times m_{d}} $ belongs to the tensor space $ \R^{(m_{1} \times n_{1}) \times \dots \times (m_{d} \times n_{d})} $. The multiplication of $\tensor{F}$ with a tensor $ \tensor{A} \in \R^{n_{1} \times \dots \times n_{d}} $ is defined as follows:
\begin{align*}
  (\tensor{F}\tensor{A})(j_{1}, \dots, j_{d}) = \sum_{i_{1}, \dots, i_{d}}\tensor{F}(j_{1}, i_{1}, \dots, j_{d}, i_{d})\tensor{A}(i_{1}, \dots, i_{d}),
\end{align*}
for all $ j_{k} = 1, \dots, m_{k}, i_{k} = 1, \dots, n_{k}$, and  $ k = 1, \dots, d $. Typically, the tensor operator $\tensor{F}$ can be represented as a tensor train operator, i.e., $ \tensor{F} = \ttcore{F}{1} \times^{1} \ttcore{F}{2} \times^{1} \dots \times^{1} \ttcore{F}{d} $, where $ \ttcore{F}{k} \in \R^{r_{k - 1} \times m_{k} \times n_{k} \times r_{k}} $ for $ k = 1, \dots, d $, with $ r_{0} = r_{d} = 1 $. We call these tensor operators in TT representation \emph{TT matrices} \cite{oseledets2011tensor}. In this case, the multiplication of $\tensor{F}$ with a TT tensor $ \tensor{A} 
$ can be decomposed into the corresponding operations on TT cores\cite{oseledets2011tensor}.


Let us consider the Hadamard product $ \tensor{A}:=\tensor{Y} \odot \tensor{Z} $, where TT tensors $\tensor{Y},\tensor{Z}\in\R^{n_1\times\cdots\times n_d}$ with TT ranks $\{r_k\}_{k=0}^d$ and $\{s_k\}_{k=0}^d$, respectively. The representation of $ \tensor{A} $ by element is 
\begin{align*}
  \tensor{A}(i_{1}, \dots, i_{d}) & = \tensor{Y}(i_{1}, \dots, i_{d})\cdot \tensor{Z}(i_{1}, \dots, i_{d}) \\
  & = (\ttcore{Y}{1}(i_{1})\cdots\ttcore{Y}{d}(i_{d}))\otimes(\ttcore{Z}{1}(i_{1})\cdots\ttcore{Z}{d}(i_{d})) \\
  & = (\ttcore{Y}{1}(i_{1}) \otimes \ttcore{Z}{1}(i_{1}))\cdots(\ttcore{Y}{d}(i_{d}) \otimes \ttcore{Z}{d}(i_{d})).
\end{align*}
Here we change the multiplication of two real number from elementwise to the Kronecker product. Due to the property $ (\mat{A}\mat{B}) \otimes (\mat{C}\mat{D}) = (\mat{A}\otimes \mat{C})(\mat{B} \otimes \mat{D}) $ and the definition of the PKP, for each parenthesis,  we have:
\begin{align*}
  \tensor{A}(i_{1}, \dots, i_{d}) = \left(\ttcore{Y}{1}\pkp^{1,3}\ttcore{Z}{1}\right)(i_{1})\cdots\left(\ttcore{Y}{d}\pkp^{1,3}\ttcore{Z}{d}\right)(i_{d}).
\end{align*}
This representation is similar to  \cref{eq:tt_decom_element}, which means TT cores of $ \tensor{A} = \tensor{Y} \odot \tensor{Z} $ can be explicitly computed by the PKP of TT cores of $ \tensor{Y} $ and $ \tensor{Z} $: $ \ttcore{A}{k} = \ttcore{Y}{k}\pkp^{1, 3}\ttcore{Z}{k} $, for all $ k = 1, \dots, d $. However, the TT ranks of $\tensor{A}$ will increase to $\{r_ks_k\}_{k=0}^d$, resulting in a significant increase of both computational complexity and memory requirements for storing $\tensor{A}$. 
Later in this paper, we denote   
$$n = \max_{1\leq k\leq d} n_{k},~r = \max_{0\leq k\leq d} r_{k},~s = \max_{0\leq k\leq d} s_{k}, \hbox{ and } ~ \ell = \max_{0\leq k\leq d} \ell_{k}.$$
We further assume that $\max\{r, s\} \lesssim \ell \ll rs$. With these assumptions,  the computational complexity and memory requirements for explicit computing the Hadamard product $\tensor{A}$ both are ${\cal O}(d n r^2s^2)$.



\subsection{Recompression of tensor train tensor}
\label{recompressionTTsec2.3}
Nearly all basic linear algebra operations on TT tensors increase the TT ranks \cite{oseledets2011tensor}. To suppress the growth of TT ranks, several recompression algorithms were proposed. The most popular recompression algorithm for TT tensors is \alg{TT-Rounding} \cite{oseledets2011tensor}. The \alg{TT-Rounding} algorithm consists of two steps: orthogonalization and recompression (typically using SVD or truncated SVD).  During the orthogonalization step, TT cores are adjusted to satisfy left orthogonality, except for the first TT core $\ttcore{A}{1}$, using QR decomposition. Subsequently, the TT cores of the approximate tensor are constructed sequentially through the truncated SVD step. The computational procedure of \alg{TT-Rounding} is detailed in Algorithm \ref{alg:TT-Rounding}.
\begin{algorithm}[!htpb]
  \caption{TT-Rounding \cite{oseledets2011tensor}}\label{alg:TT-Rounding}
  \begin{algorithmic}[1]
    \Require A $d$th-order TT tensor $ \tensor{A} $, target TT ranks $\{\ell_{k}\}_{k = 0}^{d}$.
    \Ensure TT tensor $ \tensor{B} $ with TT ranks $\{\ell_{k}\}_{k = 0}^{d}$.
    \Function{$ \tensor{B} $}{TT-Rounding}{$ \tensor{A} $, $ \{\ell_{k}\}_{k = 0}^{d} $}
      \LComment{Right-to-left orthogonalization}
      \State Set $ \ttcore{B}{d} = \ttcore{A}{d} $.
      \For{$ k = d $ down to 2}
        \State $ [\HT{\ttcore{B}{k}}\tran, \mat{R}] = \text{QR}(\HT{\ttcore{B}{k}}\tran) $. \Comment{LQ factorization of $ \HT{\ttcore{B}{k}} $}
        \State $ \VT{\ttcore{B}{k - 1}} = \VT{\ttcore{A}{k - 1}} \mat{R}\tran $.
      \EndFor 
      \State Set $ \ttcore{B}{1} = \ttcore{A}{1} $.
      \For{$ k = 1 $ to $ d - 1 $}
        \State $ [\VT{\ttcore{B}{k}}, \mat{R}] = \alg{QR}(\VT{\ttcore{B}{k}}) $. \Comment{thin QR factorization}
        \State $ [\hat{\mat{U}}, \hat{\mat{\Sigma}}, \hat{\mat{V}}] = \alg{trucated SVD}(\mat{R}, \ell_{k + 1}) $.
        \State $ \VT{\ttcore{B}{k}} = \VT{\ttcore{B}{k}}\hat{\mat{U}} $. \Comment{the U matrix to left: $ \ttcore{B}{k} \times^{1} \hat{\mat{U}} $.}
        \State $ \HT{\ttcore{B}{k + 1}} = \hat{\mat{\Sigma}}\hat{\mat{V}}\tran \HT{\ttcore{A}{k + 1}} $. \Comment{the other to right: $ \ttcore{B}{k + 1} = \hat{\mat{\Sigma}}\hat{\mat{V}}\tran \times^{1} \ttcore{A}{k + 1} $}
      \EndFor
    \EndFunction
  \end{algorithmic}
\end{algorithm}

In the procedure of \alg{TT-Rounding}, the total computational cost of a TT tensor with rank $ r_{k} \le r $  is $ O(dnr^{3}) $. Therefore, by using \alg{TT-Rounding} to recompress the Hadamard product $ \tensor{A} = \tensor{Y} \odot \tensor{Z} $, the total computational cost amounts to $ O(dnr^{3}s^{3}) $. This could be computationally expensive, particularly for tensors $\tensor{Y}$ and $\tensor{Z}$ with relatively large ranks. 
This motivates us to develop more efficient algorithms for recompressing Hadamard products. 

The computational complexity analysis reveals that the SVD of the cores is the most computationally expensive step in \alg{TT-Rounding}. To improve efficiency, and inspired by randomized SVD, Daas \emph{et al.} \cite{al2023randomized} proposed three randomized TT recompression algorithms: Orthogonalize-then-Randomize (OrthRand), Randomize-then-Orthogonalize (RandOrth), and Two-Sided-Randomization (TwoSided). 
The core idea of these randomized TT recompression methods is to replace the truncated SVD with a more efficient randomized SVD. The primary distinction between these randomized TT recompression algorithms lies in how the sketch used in randomized SVD is generated. 
As an example, we introduce \alg{RandOrth} below.

Let us consider a matrix $ \mat{A} = \mat{B} \mat{C} $, where $ \mat{A} \in \R^{m \times n} $, $ \mat{B} \in \R^{m \times R} $, $ \mat{C} \in \R^{R \times n} $ with $ R > \rank(\mat{A}) $. 
 The goal of randomized SVD is to find an approximation of $ \mat{A}  $ with a target rank $ \ell<R $. To achieve this, we first generate a random matrix $\mat{\Omega}\in\R^{n\times \ell}$,  whose elements are independently drawn from the standard normal distribution. And then we compute the QR factorization of
\begin{align}
\label{random-qr}
 \mat{A}\mat{\Omega} = \mat{B}\mat{C}\mat{\Omega} := \mat{B}\mat{W} = \mat{Q}\mat{R},
\end{align}
where $ \mat{Q} \in \R^{m \times \ell} $ has orthogonal columns. We have $ \Range \mat{Q} \approx \Range (\mat{A}) $ and use $ 
\mat{Q}(\mat{Q}\tran \mat{A})  $ to approximate $ \mat{A} $, whose accuracy strongly depends on the singular value distribution of $ \mat{A} $ \cite{halko2011finding}. 
Inspired by randomized SVD, Daas \emph{et al.} \cite{al2023randomized} introduced the concept of a random TT tensor and used it to develop a randomized TT recompression algorithm.

\begin{definition}[Random TT tensor]\label{def:randTT-2}
  Given a set of target TT ranks $ \{\ell_{k}\}_{k = 0}^{d} $, we generate a random Gaussian TT tensor $ \tensor{R} \in \R^{n_{1} \times \dots \times n_{d}} $ such that each core tensor $ \ttcore{R}{k} \in \R^{\ell_{k - 1} \times n_{k} \times \ell_{k}} $ is filled with 
  random, independent, normally distributed entries with mean $ 0 $ and variance $ 1/(\ell_{k - 1}n_{k}\ell_{k}) $ for $ 1 \le k \le d $. If each core is filled with random, independent, uniformly distributed entries from $ [0, 1] $, we call it a random uniform TT tensor.
\end{definition}

Let us rewrite the partial contracted product $ \ttcore{A}{k:d} $ as two tensors' multiplication: $ \ttcore{A}{k} \times^{1} \ttcore{A}{(k+1):d} $. Similar to \cref{random-qr}, we multiply a random tensor  $\tensor{R}^{(k)}\in\R^{n_{k+1} \times n_{k+2} \times \dots \times n_{d} \times \ell_{k}}$ on the right-hand side to generate a sketch, i.e., 
\begin{align*}
 \mat{W}^{(k)} = \ttcore{A}{(k+1):d} \times^{1, 2, \dots, d - k} \tensor{R}^{(k)} \in \R^{r_{k}s_{k} \times \ell_{k}}.
\end{align*}
If $\tensor{R}^{(k)}$ is a sub-tensor of a random Gaussian TT tensor $\tensor{R}$, i.e.,
\begin{align*}
    \tensor{R}^{(k)} = \ttcore{R}{k + 1} \times^{1} \ttcore{R}{k + 2} \times^{1} \dots \times^{1} \ttcore{R}{d},
\end{align*}
then the corresponding randomized TT recompression algorithm is RandOrth proposed in \cite{al2023randomized}, see Algorithm \ref{alg:RandOrth}.
\begin{algorithm}[!htpb]
  \caption{RandOrth: Randomize-then-Orthogonalize \cite{al2023randomized}}\label{alg:RandOrth}
  \begin{algorithmic}[1]
    \Require A TT tensor $ \tensor{A} $ with ranks $ \{r_{k}\}_{k = 0}^{d} $, target TT ranks $ \{\ell_{k}\}_{k = 0}^{d} $ 
    \Ensure A TT tensor $ \tensor{B} $ with ranks $ \{\ell_{k}\}_{k = 0}^{d} $
    \Function{$ \tensor{B} $}{RandOrth}{$ \tensor{A} $, $ \{\ell_{k}\}_{k = 0}^{d} $}
      \State Select a random Gaussian TT tensor $ \tensor{R} $ with target TT ranks $ \{\ell_{k}\}_{k = 0}^{d} $
      \State $ \mat{W}^{(d - 1)} = \HT{\ttcore{A}{d}}\left(\HT{\ttcore{R}{d}}\right)\tran $
      \For{$ k = d - 1 $ down to $ 2 $}
        \State $ \VT{\ttcore{Z}{k}} = \VT{\ttcore{A}{k}}\mat{W}^{(k)} $
        \State $ \mat{W}^{(k - 1)} = \HT{\ttcore{Z}{k}}\left(\HT{\ttcore{R}{k}}\right)\tran $ 
      \EndFor
      \State $ \ttcore{B}{1} = \ttcore{A}{1} $
      \For{$ k = 1 $ to $ d - 1 $} 
        \State $ \mat{B}^{(k)} = \VT{\ttcore{B}{k}} $ 
        \State $ [\VT{\ttcore{B}{k}}, \sim] = \text{QR}(\mat{B}^{(k)}\mat{W}^{(k)}) $
        \State $ \mat{M}^{(k)} = \left(\VT{\ttcore{B}{k}}\right)\tran \mat{B}^{(k)} $
        \State $ \HT{\ttcore{B}{k + 1}} = \mat{M}^{(k)}\HT{\ttcore{A}{k + 1}} $
      \EndFor
    \EndFunction
  \end{algorithmic}
\end{algorithm}

According to the numerical simulations reported in \cite{al2023randomized}, the \alg{RandOrth} algorithm, along with \alg{OrthRand} and \alg{TwoSided}, outperforms traditional TT-rounding. These methods are effective for recompressing both individual TT tensors and the summations of multiple TT tensors.
However, when applied to recompress the Hadamard product of multiple TT tensors, these methods still face challenges in terms of computational complexity and memory usage, particularly for TT tensors with relatively large ranks. To address these limitations, we propose a new algorithm in the next section, which is inspired by these randomized TT recompression methods and designed to enhance the efficiency and performance of recompressing the Hadamard product of multiple TT tensors.

\section{Hadamard avoiding TT  recompression}
\label{sec:hatt}


Our goal is to efficiently round the Hadamard product $ \tensor{A}=\tensor{Y} \odot \tensor{Z} $ into a TT tensor with ranks $ \{\ell_{k}\}_{k=0}^d $. As discussed in \autoref{recompressionTTsec2.3}, existing recompression algorithms require explicitly computing the Hadamard product $ \tensor{A}=\tensor{Y} \odot \tensor{Z} $. This computation entails both computational complexity and memory requirements of ${\cal O}(d n r^2s^2)$, which becomes prohibitively expensive for large values of $r$ and $s$.  
Inspired by RandOrth \cite{al2023randomized} (see \cref{alg:RandOrth}), we propose the HaTT algorithm recompresses $\tensor{A}$ using the partial contraction of $\tensor{A}$ and a random Gaussian TT tensor $\tensor{R}$. The HaTT algorithm is more efficient as it utilizes the PKP structure in TT cores of $\tensor{A}$ and avoids the explicit computation of TT tensor $\tensor{A}$. 
The HaTT algorithm starts by generating sketch matrices through the partial contraction of the Hadamard product with a random tensor. It then performs QR factorization on these sketch matrices for each core to obtain the resulting low-rank TT tensor.
To achieve this, we first introduce the partial contraction process for the Hadamard product to generate sketched matrices in \autoref{subsec:HPCRL}. Next, we present the HaTT algorithm by detailing the QR factorization of the sketched matrices in \autoref{QRforeachcore}.

\subsection{Right-to-left partial contraction for Hadamard product}\label{subsec:HPCRL}
We begin by introducing the partial contraction of TT tensors $\tensor{A}:=\tensor{Y}\odot \tensor{Z}$ and $\tensor{R}$. According to \cite{al2023randomized}, for $k=2, 3, \dots, d$, the partial contraction matrices are defined by:
\begin{align}\label{eq:partial-contraction}
  \mat{W}^{(k-1)} & = \HT{\ttcore{A}{k:d}}\left(\HT{\ttcore{R}{k:d}}\right)\tran\in \R^{r_{k-1}s_{k-1}\times \ell_{k-1}}, 
\end{align} 
and satisfy the following recursion formula
\begin{equation}
\begin{aligned}
\label{rtlpc}
  \VT{\ttcore{B}{k}} & = \VT{\ttcore{A}{k}}\mat{W}^{(k)}\in \R^{r_{k-1}s_{k-1}n_{k}\times \ell_{k}}\\ 
  \mat{W}^{(k-1)} & = \HT{\ttcore{B}{k}}\left(\HT{\ttcore{R}{k}}\right)\tran. 
\end{aligned}
\end{equation}
Here, $\tensor{B}$ is a temporary TT tensor with compatible dimensions and ranks. The recursion formula \cref{rtlpc} can be rewritten in matrix form 
\begin{align*}
  \mat{W}^{(k-1)}:= \sum_{i_{k} = 1}^{n_{k}}\ttcore{A}{k}(i_{k}) \mat{W}^{(k)} \left(\ttcore{R}{k}(i_{k})\right)\tran.
\end{align*}
The process of computing the matrices $\{\mat{W}^{(k-1)}\}_{k=2}^d$ according to \cref{rtlpc} is called {right-to-left partial contraction (PartialContractionRL) \cite{al2023randomized}, which is displayed in \cref{fig:randorth-process}. 
If the TT tensor $\tensor{A}$ is a result of Hadamard product of $ 
\tensor{Y} $ and $ \tensor{Z} $, the total computational cost for the PartialContractionRL is $O(dn(rs\ell^2+r^2s^2\ell))$  \cite{al2023randomized}.

\begin{figure}
  \centering
  \begin{tikzpicture}[baseline]
  \Vertex[x=-.5, y=0.75, size=0, Pseudo=True]{Ak-1}
  \Vertex[x=-.5, y=-0.75, size=0, Pseudo=True]{Rk-1}
  \Vertex[x=1, y=0.75, label={$ \ttcore{A}{k} $}, position=above, distance=-1mm]{Ak}
  \Vertex[x=1, y=-0.75, label={$ \ttcore{R}{k} $}, position=below, distance=-1mm]{Rk}
  \Vertex[x=2, y=0, label={$ \mat{W}^{(k)} $}, position=right]{Wk}
  \Edge(Ak-1)(Ak)
  \node[below=2pt, color={black!75}] at (-0.1, 0.75) {\tiny $ r_{k - 1}s_{k - 1} $};
  \Edge[label=$ \ell_{k - 1} $, position={above}](Rk-1)(Rk)
  \Edge[label=$ r_{k}s_{k} $, position={above right}](Ak)(Wk)
  \Edge[label=$ \ell_{k} $, position={below right}](Rk)(Wk)
  \Edge[label=$ n_{k} $, position=left](Ak)(Rk)
\end{tikzpicture} $ \longrightarrow $ 
\begin{tikzpicture}[baseline]
  \Vertex[x=-.5, y=0.75, size=0, Pseudo=True]{Ak-1}
  \Vertex[x=-.5, y=-0.75, size=0, Pseudo=True]{Rk-1}
  \Vertex[x=1, y=0.75, label={$ \ttcore{B}{k} $}, position=above, distance=-1mm]{Zk}
  \Vertex[x=1, y=-0.75, label={$ \ttcore{R}{k} $}, position=below, distance=-1mm]{Rk}
  \Edge(Zk)(Ak-1)
  \node[below=2pt, color={black!75}] at (-0.1, 0.75) {\tiny $ r_{k - 1}s_{k - 1} $};
  \Edge[label=$ \ell_{k-1} $, position={above}](Rk)(Rk-1)
  \Edge[label=$ n_{k} $, bend=-30, position=left](Zk)(Rk)
  \Edge[label=$ \ell_{k} $, bend=30, position=right](Zk)(Rk)
\end{tikzpicture} $ \longrightarrow $ 
\begin{tikzpicture}[baseline]
  \Vertex[x=-.5, y=0.75, size=0, Pseudo=True]{Ak-1}
  \Vertex[x=-.5, y=-0.75, size=0, Pseudo=True]{Rk-1}
  \Vertex[x=1, y=0, label={$ \mat{W}^{(k - 1)} $}, , position=right]{Wk-1}
  \Edge(Wk-1)(Ak-1)
  \Text[x=-.3, y=.75, fontsize=\tiny, position={above right}, rotation={-25}, color=black!75]{$ r_{k - 1}s_{k - 1} $}
  \Edge[label=$ \ell_{k - 1} $, position={below right=.5pt}](Wk-1)(Rk-1)
\end{tikzpicture} 
  \caption{The process of right-to-left partial contraction}\label{fig:randorth-process}
\end{figure}

To reduce the computational cost of the PartialContractionRL for the case $\tensor{A}:=\tensor{Y}\odot\tensor{Z}$, we introduce a new approach with the Hadamard avoiding technique in this subsection. This method utilizes the PKP operation and property \cref{eq:kronecker_product_times_vec} to reformulate the recursion formula \cref{rtlpc}, enabling the partial contraction to be computed without explicitly constructing the TT core tensor $\ttcore{A}{k}$. The reformulated recursion formula \cref{rtlpc} is derived through the following two steps. 

\vspace{2mm}

\begin{enumerate}[label={\bf{Step \arabic*:}}, leftmargin=0mm, itemindent=*, widest=2]
  \item Represent the matrix $\mat{W}^{(k)}\in\R ^{r_{k}s_{k}\times \ell_k}$ as the summation of multiple rank-1 matrices, i.e. $ \mat{W}^{(k)} = \sum_{\gamma=1}^{\ell_k}\sigma^{(k)}_{\gamma}\mat{U}^{(k)}(:,{\gamma})\left(\mat{V}^{(k)}(:,{\gamma})\right)\tran:=\mat{U}^{(k)}\mat{S}^{(k)}\left(\mat{V}^{(k)}\right)\tran $, where $ \mat{U}^{(k)} \in \R^{r_{k}s_{k} \times \ell_{k}} $, $ \mat{V}^{(k)} \in \R^{\ell_{k} \times \ell_{k}} $, and $ \mat{S}^{(k)} := \diag\{\sigma_{1}^{(k)}, \sigma_{2}^{(k)}, \dots,\sigma_{\ell_k}^{(k)}\} \in \R^{\ell_{k} \times \ell_{k}} $. 
  Then we reshape vector $ \mat{U}^{(k)}(:, \gamma) \in \R^{r_{k}s_{k}} $ to a matrix $ \mat{U}^{(k)}_{\gamma}\in \R^{s_{k} \times r_{k}}  $.
  \item The partial contraction matrix $ \mat{W}^{(k-1)} $ is calculated by:
  \begin{small}
  \begin{align}
    \mat{W}^{(k-1)} & = \sum_{i_{k} = 1}^{n_{k}}\ttcore{A}{k}(i_{k}) \mat{W}^{(k)} \left(\ttcore{R}{k}(i_{k})\right)\tran \notag\\
    & = \sum_{i_{k} = 1}^{n_{k}}(\ttcore{Y}{k}(i_{k}) \otimes \ttcore{Z}{k}(i_{k}))\mat{W}^{(k)} \left(\ttcore{R}{k}(i_{k})\right)\tran \notag\\
    & = \sum_{i_{k} = 1}^{n_{k}}(\ttcore{Y}{k}(i_{k}) \otimes \ttcore{Z}{k}(i_{k})) \mat{U}^{(k)}\mat{S}^{(k)}\left(\mat{V}^{(k)}\right)\tran \left(\ttcore{R}{k}(i_{k})\right)\tran \label{eq:partial contraction with svd}\\
    & = \sum_{i_{k} = 1}^{n_{k}} \left(\sum_{\gamma = 1}^{\ell_{k}}\sigma_{\gamma}^{(k)}\big(\ttcore{Y}{k}(i_{k}) \otimes \ttcore{Z}{k}(i_{k})\big)\mat{U}^{(k)}(:, \gamma) \left(\mat{V}^{(k)}(:, \gamma)\right)\tran \left(\ttcore{R}{k}(i_{k})\right)\tran \right)\notag\\
    & = \sum_{i_{k} = 1}^{n_{k}}\sum_{\gamma = 1}^{\ell_{k}} \sigma_{\gamma}^{(k)} \alg{vec}\big(\ttcore{Z}{k}(i_{k}) \mat{U}^{(k)}_{\gamma} (\ttcore{Y}{k}(i_{k}))\tran\big) \big(\ttcore{R}{k}(i_{k})\mat{V}^{(k)}(:, \gamma)\big)\tran, \notag
  \end{align}
  \end{small}
  where the last equation holds due to property \cref{eq:kronecker_product_times_vec}.
\end{enumerate}

The process of computing the matrices $\{\mat{W}^{(k-1)}\}_{k=2}^d$ according to Steps 1 and 2 is referred to as PartialContractionRL for Hadamard product (HPCRL). 
According to \cref{eq:partial contraction with svd}, the partial contraction $ \mat{W}^{(k-1)} $ is derived from matrix-matrix multiplication of slices $\ttcore{Z}{k}(i_{k})$ and $\ttcore{Y}{k}(i_{k})$, bypassing the need to explicitly compute $\ttcore{A}{k}(i_{k}) $. 
To streamline this computation, we introduce two matrices, $ \mat{W}_{L} \in \R^{r_{k - 1}s_{k - 1} \times n_{k}\ell_{k}} $ and $ \mat{W}_{R} \in \R^{\ell_{k - 1} \times n_{k}\ell_{k}} $, defined as:
\begin{equation}
\label{WLWR}
\left\{\begin{aligned}
  \mat{W}_{L}(:, \overline{i_{k}\gamma}) &= \alg{vec}\big(\ttcore{Z}{k}(i_{k}) \mat{U}^{(k)}_{\gamma} (\ttcore{Y}{k}(i_{k}))\tran\big),\\  \mat{W}_{R}(:, \overline{i_{k}\gamma}) &= \ttcore{R}{k}(i_{k})\mat{V}^{(k)}(:, \gamma),
  \end{aligned}
  \right.
  \end{equation}
where $\overline{i_{k}\gamma}$ denotes the column indexing. 
With these definitions, \cref{eq:partial contraction with svd} can be rewritten as:
\begin{equation}\label{LISRmatrix}
\mat{W}^{(k-1)} = \mat{W}_{L}(\mat{I}_{n_{k}} \otimes \mat{S}^{(k)})\mat{W}_{R}\tran.
\end{equation}
This matrix formulation is straightforward to implement in MATLAB and takes advantage of efficient matrix-matrix multiplication, significantly enhancing computational performance.

In Step 1 of HPCRL, the matrix $\mat{W}^{(k)}\in\R ^{r_{k}s_{k}\times \ell_k}$ must be expressed as a summation of multiple rank-1 matrices. There are many ways to achieve this representation, and the choice depends on balancing the computational cost of the representation against the sparsity of the diagonal matrix $\mat{S}^{(k)}$. 
On one hand, the computational cost of Step 2 is directly proportional to the number of non-zero elements in $\mat{S}^{(k)}$. Therefore, a sparser $\mat{S}^{(k)}$ leads to greater efficiency in the proposed approach. On the other hand, the computational cost of constructing the representation increases as the number of non-zero elements in $\mat{S}^{(k)}$ decreases. 
In this paper, we present two methods to achieve this representation, providing options to balance efficiency and computational complexity. 


The first method directly represents the matrix $\mat{W}^{(k)}$ as
\begin{align}
\label{HPCRL--1}
\mat{W}^{(k)}=\sum_{\gamma=1}^{\ell_k}{\mat{U}^{(k)}(:,{\gamma})}\left(\mat{V}^{(k)}(:,{\gamma})\right)\tran,
\end{align}
where $ \mat{U}^{(k)}(:, \gamma) $ and $\mat{V}^{(k)}(:,{\gamma})$ are the $\gamma$th column of $\mat{W}^{(k)}$ and $\mat{I}_{\ell_k}$, respectively. This representation is almost free, but the sparsity of $\mat{S}^{(k)} = \mat{I}_{\ell_{k}}$ is the worst. {\color{black} Since $ \mat{S} = \mat{I}_{n\ell} $, we can omit the Kronecker product $ \mat{I}_{n_{k}} \otimes \mat{S} $ in \eqref{LISRmatrix}.} Compared to the standard partial contraction, the computational cost of HPCRL is reduced from $O(dn\cdot rs\ell(rs + \ell))$ to $O(dn\cdot rs\ell(r + s + \ell))$.
    
The second method is (truncated) SVD, represented as: 
    \begin{align}
    \label{HPCRL--2}
        \mat{W}^{(k)} = \sum_{\gamma=1}^{R}\sigma^{(k)}_{\gamma} \mat{U}^{(k)}(:,{\gamma})\left(\mat{V}^{(k)}(:,{\gamma})\right)\tran=\color{black}\sum_{\gamma=1}^{R} \mat{U}^{(k)}(:,{\gamma})\left(\sigma^{(k)}_{\gamma}\mat{V}^{(k)}(:,{\gamma})\right)\tran,
    \end{align}
    where $R$ is either $\rank(\mat{W}^{(k)})$ or the target rank of $\mat{W}^{(k)}$ in the case of truncated SVD. 
 The computational cost to perform SVD is $O(rs\ell^{2})$. {\color{black} By multiplying $ \sigma^{(k)}_{\gamma} $ with $ \mat{V}^{(k)}(:,{\gamma}) $, the Kronecker product $ \mat{I}_{n_{k}} \otimes \mat{S} $ in \eqref{LISRmatrix} can also be omitted.}
Given that $R\leq \ell_k$ and the use of SVD in \eqref{HPCRL--2}, the summations in the last two equations of \eqref{eq:partial contraction with svd} for $\gamma$ are limited to $R$. Consequently, the  resulting $\mat{W}_L$ is a matrix of size $\mathbb{R}^{r_{k-1}s_{k-1}\times n_kR}$, while $\mat{W}_R$ is of size $\mathbb{R}^{\ell_{k-1} \times n_kR}$. Compared to the standard partial contraction, the computational cost of HPCRL decreases from $O(dnrs\ell(rs + \ell))$ to  $O(drs  [nR(r + s + \ell)+\ell^2])$.

The HPCRL process using the direct representation outlined in \cref{HPCRL--1} is illustrated in \cref{fig:HPCRL-No-SVD} and referred to as HPCRL-1. The HPCRL process employing the (truncated) SVD method described in \cref{HPCRL--2} is depicted in \cref{fig:HPCRL-process} and designated as HPCRL-2. A comprehensive summary of the HPCRL process is provided in \cref{alg:HPCRL-main}. 
{\color{black} The advantages and disadvantages of HPCRL-1 and HPCRL-2 depend on the rank of $\mat{W}^{(k)}$ and the target rank $\ell$. 
To elucidate this relationship, we now present a theoretical analysis.
Specifically, the rank of $\mat{W}^{(k)} $ is bounded above by the $ k $-th separation rank of the TT tensor $ \tensor{A} $, i.e., the rank of its unfolding matrix $ \unfold{A}{k} $. This result is formally stated in the following theorem.
\begin{theorem}
\label{theorem:upboundofW}
  Let $ \tensor{A} \in \R^{n_{1} \times n_{2} \times \cdots \times n_{d}} $ be a TT tensor with separation ranks $ \{r^{\text{sep}}_{0}, r^{\text{sep}}_{1}, \dots, r^{\text{sep}}_{d}\} $. Suppose that we generate a set of sketch matrices $ \{\mat{W}^{(k)}\}_{k = 1}^{d - 1} $ of $ \tensor{A} $ using a random TT tensor $ \tensor{R} $ through equation \cref{eq:partial-contraction}. Then for any $ k = 1, 2, \dots, d - 1 $, it holds that $ \rank \mat{W}^{(k)} \leq r^{\text{sep}}_{k} $.
\end{theorem}
\begin{proof}
  By the definition of the separation rank \cite{holtz2012alternating}, the unfolding matrix $ \unfold{A}{k} $ admits the decomposition 
  \begin{align*}
    \unfold{A}{k} = \VT{\ttcore{A}{1:k}} \HT{\ttcore{A}{k + 1:d}},
  \end{align*}
  where $ \rank \unfold{A}{k} = \rank \VT{\ttcore{A}{1:k}} = \rank \HT{\ttcore{A}{k + 1:d}} = r^{\text{sep}}_{k} $. 
  According to equation~\cref{eq:partial-contraction}, we have
  \begin{align*}
    \rank \mat{W}^{(k - 1)} &\le \min\left\{\rank \HT{\ttcore{A}{k:d}}, \rank \HT{\ttcore{R}{k:d}}\right\} \\ 
    &\le \rank \HT{\ttcore{A}{k:d}} =  r^{\text{sep}}_{k-1},
  \end{align*}
  which completes the proof.
\end{proof}


When performing TT recompression, our goal is to compress a TT tensor with relatively high separation ranks into a low-rank TT representation.
Following the definition for matrices in \cite{beckermann2017singular}, we define \emph{quasi-numerical ranks} of a TT tensor as follows. 
\begin{definition}[Absolute quasi-numerical ranks]
  Let $ \tensor{A} \in \R^{n_{1} \times n_{2} \times \cdots \times n_{d}} $ be a TT tensor and $ \varepsilon > 0 $ be a given absolute error tolerance. Let $ \tensor{B} $, with separation ranks $\{r_{0}^{\varepsilon}, r_{1}^{\varepsilon}, \dots, r_{d}^{\varepsilon}\}$, denote the result of applying the TT-SVD algorithm~\cite{oseledets2011tensor} to $ \tensor{A} $, such that $ \|\tensor{A} - \tensor{B}\|_{F} \leq \varepsilon  $. Then $\{r_{0}^{\varepsilon}, r_{1}^{\varepsilon}, \dots, r_{d}^{\varepsilon}\}$ is called the \emph{absolute quasi-numerical ranks} of $ \tensor{A} $ with respect to $ \varepsilon $, denoted by $ \rank_{a,\varepsilon}(\tensor{A}) = \{r_{0}^{\varepsilon}, r_{1}^{\varepsilon}, \dots, r_{d}^{\varepsilon}\} $.
\end{definition}    

\begin{definition}[Relative quasi-numerical ranks]
  Let $ \tensor{A} \in \R^{n_{1} \times n_{2} \times \cdots \times n_{d}} $ be a TT tensor and $ \varepsilon > 0 $ be a given relative error tolerance. Let $ \tensor{B} $, with separation ranks $\{r_{0}^{\varepsilon}, r_{1}^{\varepsilon}, \dots, r_{d}^{\varepsilon}\}$, denote the result of applying the TT-SVD algorithm~\cite{oseledets2011tensor} to $ \tensor{A} $, such that $ \|\tensor{A} - \tensor{B}\|_{F} \leq \varepsilon \|\tensor{A}\|_F  $. Then $\{r_{0}^{\varepsilon}, r_{1}^{\varepsilon}, \dots, r_{d}^{\varepsilon}\}$ is called the \emph{relative quasi-numerical ranks} of $ \tensor{A} $ with respect to $ \varepsilon $, denoted by $ \rank_{r,\varepsilon}(\tensor{A}) = \{r_{0}^{\varepsilon}, r_{1}^{\varepsilon}, \dots, r_{d}^{\varepsilon}\} $.
\end{definition}   
\begin{remark}
    For matrices, the numerical rank can be determined using SVD. However, for TT tensors, the numerical ranks, similar to border ranks with given tolerance \cite{bini2007role}, are generally difficult to compute. In this work, we define the numerical TT ranks obtained by the TT-SVD algorithm with a prescribed tolerance as quasi-numerical ranks, to distinguish them from border ranks. In the following, we use the term ``numerical rank” interchangeably with ``quasi-numerical rank” for simplicity.
\end{remark}

\begin{lemma}
\label{lemma:lowerboundofA-B}
  Let $ \tensor{A} $ be a left orthogonal TT tensor with absolute numerical ranks $ \rank_{a,\varepsilon}(\tensor{A}) $ $= \{r_{0}^{\varepsilon}, r_{1}^{\varepsilon}, \dots, r_{d}^{\varepsilon}\} $. Let $ \tensor{B} $ denote the result of applying the TT-SVD algorithm~\cite{oseledets2011tensor} to $ \tensor{A} $, such that $ \|\tensor{A} - \tensor{B}\|_{F} \leq \varepsilon  $. Then we have 
    \begin{align}
    \label{lemma:eq-2}
\left\|\tensor{A}-\tensor{B}\right\|_F^2\geq \sum\limits_{k=1}^{d-1} 
    \sigma^2_{r_k^\varepsilon+1}(\HT{\ttcore{A}{k + 1:d}}),
  \end{align}
  where $\sigma_{r_k^\varepsilon+1}(\HT{\ttcore{A}{k + 1:d}})$ denotes the $(r_k^\varepsilon+1)$-th singular value of $\HT{\ttcore{A}{k + 1:d}}$. Therefore, the absolute $\varepsilon$-numerical rank of $\HT{\ttcore{A}{k + 1:d}}$ does not exceed $r_k^\varepsilon$.
\end{lemma}
\begin{proof}
 Following the proof of Theorem 2.2 in \cite{oseledets2011tensor}, the argument proceeds by induction. For $d=2$, the statement follows directly from the properties of the SVD. To simplify the notation, denote $\mat{V}_{k}=\VT{\ttcore{A}{1:k}}$ and $\mat{H}_{k+1}=\HT{\ttcore{A}{k+1:d}}$.  For an arbitrary $d>2$, the first unfolding matrix $\unfold{A}{1}$ can be decomposed as
 \begin{align*}
    \unfold{A}{1}=\mat{V}_{1}\mat{H}_{2}=\mat{V}_{1}\mat{\tilde{V}}_{1}\bm{\Sigma}_{r_1^\varepsilon}\mat{\bar{U}}_{1}+\unfold{E}1:=\mat{\bar{V}}_{1}\mat B_1+\unfold{E}1,
  \end{align*}
  where $\|\unfold{E}1\|_F\geq \sigma_{r_1^\varepsilon+1}(\mat{H}_{2})$, and $\mat{\bar{V}}_{1}:=\mat{V}_{1}\mat{\tilde{V}}_{1}$ is of size $n_1\times r_1^\varepsilon$. The matrix $\mat{B}_1$ corresponds naturally to a $(d-1)$-dimensional tensor $\tensor{B}_1$, which is further decomposed in the TT-SVD algorithm. Hence, $\mat{B}_1$ is approximated by another matrix $\mat{\hat{B}}_1$. Due to the orthogonality of the columns of $\mat{\bar{V}}{1}$, it follows that $\mat{\bar{V}}_{1}\tran \unfold{E}{1}=0$, and therefore 
 \begin{equation}
 \begin{aligned}\label{lemma:eq-1}
 \|\tensor{A} - \tensor{B}\|^2_{F} &= \|\unfold{A}{1}-\mat{\bar{V}}_1\mat{\hat{B}}_1\|_F^2=\left\|\unfold{A}{1}-\mat{\bar{V}}_1\left(\mat{\hat{B}}_1+\mat{B}_1-\mat{B}_1\right)\right\|_F^2\\&
 =\left\|\unfold{A}{1}-\mat{\bar{V}}_1\mat{B}_1\right\|_F^2+\left\|\mat{\bar{V}}_1(\mat{B}_1-\mat{\hat{B}}_1)\right\|_F^2\\
 &\geq \sigma^2_{r_1^\varepsilon+1}(\mat{H}_{2}) + \left\|\mat{B}_1-\mat{\hat{B}}_1\right\|_F^2.
 \end{aligned}
 \end{equation}
Since the columns of $\mat{\bar{V}}_1$ are orthonormal, the unfolding matrix $\unfold{B}{k}$ of the $(d-1)$-dimensional tensor $\tensor{B}_1$ has a similar structure to that of  $\tensor{A}$, namely 
$$
\unfold{B}{k}=\VT{\ttcore{B}{1:k}} \HT{\ttcore{B}{k + 1:d-1}},
$$
where the singular values of $\HT{\ttcore{B}{k + 1:d-1}}$ coincide with those of $\HT{\ttcore{B}{k :d}}$. Proceeding by induction, we obtain
$$
\left\|\mat{B}_1-\mat{\hat{B}}_1\right\|_F^2\geq \sum\limits_{k=2}^{d-1} \sigma^2_{r_k^\varepsilon+1}(\mat{H}_{k+1}),
$$
and combining this with \eqref{lemma:eq-1} completes the proof of \eqref{lemma:eq-2}.

Now, suppose that the $\varepsilon$-numerical rank of $\HT{\ttcore{A}{k+1:d}}$ exceeds $r_k^\varepsilon$, implying that
$$
\sigma_{r_k^\varepsilon+1}(\HT{\ttcore{A}{k + 1:d}})> \varepsilon.
$$
Substituting this into \eqref{lemma:eq-2} yields
$$
\|\tensor{A}-\tensor{B}\|> \varepsilon,
$$
which contradicts the condition $\|\tensor{A}-\tensor{B}\|\leq \varepsilon$. Hence, the proof is complete.
\end{proof}

\begin{lemma}[Theorem 3.3.16 in \cite{horn1994topics}]
 \label{multiplicative Weyl inequality}
  The multiplicative Weyl inequality for the singular values of matrices $H_1\in \mathbb{C}^{n\times m}$ and $H_2\in \mathbb{C}^{m\times n}$ are:
  \begin{align}
    \sigma_{i+j-1}(H_1H_2)\leq \sigma_{i}(H_1)\sigma_j(H_2),
\end{align}
where $1\leq i,j\leq \min\{n,\,m\} $, and $ i+j-1\leq \min\{n,\,m\}$. 

\end{lemma}

\begin{theorem}
\label{theorem3.5numericallowrankofH}
  Let $ \tensor{A} =\tensor{Y}\odot \tensor{Z}$ with $\tensor{Y}$ and $ \tensor{Z}$ being left orthogonal TT tensors. The $\varepsilon$-numerical ranks for $\tensor{A}$ are $ \rank_{a,\varepsilon}(\tensor{A}) $ $= \{r_{0}^{\varepsilon}, r_{1}^{\varepsilon}, \dots, r_{d}^{\varepsilon}\} $. Let $ \tensor{B} $ denote the result of applying the TT-SVD algorithm~\cite{oseledets2011tensor} to $ \tensor{A} $, such that $ \|\tensor{A} - \tensor{B}\|_{F} \leq \varepsilon  $. Assume that the matrices $\VT{\ttcore{A}{1:k}}$ all are column-full rank. Then we have 
    \begin{align}
\label{theorem3.5:eq-0}
\left\|\tensor{A}-\tensor{B}\right\|_F^2\geq \sum\limits_{k=1}^{d-1} 
    \frac{\sigma^2_{r_k^\varepsilon+j}(\HT{\ttcore{A}{k + 1:d}})}{\sigma^2_{j}\left(\VT{\ttcore{A}{1:k}}^+\right)},
  \end{align}
  where $\sigma_{r_k^\varepsilon+1}(\HT{\ttcore{A}{k + 1:d}})$ denotes the $r_k^\varepsilon+1$ singular value of $\HT{\ttcore{A}{k + 1:d}}$ and $\VT{\ttcore{A}{1:k}}^+$ represents the Moore-Penrose inverse of $\VT{\ttcore{A}{1:k}}$. Then, the absolute $\frac\varepsilon {\sigma_{j}\left(\VT{\ttcore{A}{1:k}}^+\right)}$ numerical rank of $\HT{\ttcore{A}{k + 1:d}}$ does not exceed $r_k^\varepsilon+j-1$.
\end{theorem}

\begin{proof}
 It is straightforward from Theorem 2.2 in \cite{oseledets2011tensor} or \Cref{lemma:lowerboundofA-B} that
 \begin{align}
 \label{theorem3.5:eq-1}
\left\|\tensor{A}-\tensor{B}\right\|_F^2\geq \sum\limits_{k=1}^{d-1} 
    \sigma^2_{r_k^\varepsilon+1}\left(\VT{\ttcore{A}{1:k}}\HT{\ttcore{A}{k + 1:d}}\right).
  \end{align}
   According to \Cref{multiplicative Weyl inequality}, we have
   \begin{equation}
 \label{theorem3.5:eq-2}
    \begin{aligned}
    \sigma_{j}\left(\VT{\ttcore{A}{1:k}}^+\right)&\sigma_{r_k^\varepsilon+1}\left(\VT{\ttcore{A}{1:k}}\HT{\ttcore{A}{k + 1:d}}\right)\\&\geq \sigma_{j+r_k^\varepsilon}\left(\VT{\ttcore{A}{1:k}}^+\VT{\ttcore{A}{1:k}}\HT{\ttcore{A}{k + 1:d}}\right)\\&= \sigma_{j+r_k^\varepsilon}\left(\HT{\ttcore{A}{k + 1:d}}\right).
  \end{aligned}
   \end{equation}
 Combining \eqref{theorem3.5:eq-1} and \eqref{theorem3.5:eq-2}, we obtain \eqref{theorem3.5:eq-0}. The proof can then be completed in a manner analogous to that of \Cref{lemma:lowerboundofA-B}. 
\end{proof}


If $\tensor{A}$ has relatively low numerical ranks, \Cref{theorem3.5numericallowrankofH} suggests that the matrix $\HT{\ttcore{A}{k + 1:d}}$ may also be numerically low rank, weighted by the singular values of $\VT{\ttcore{A}{1:k}}^+$. Since $\VT{\ttcore{A}{1:k}}$ is obtained from the matricization of the Hadamard product of two left orthogonal tensors, its singular values decay more slowly than those of $\HT{\ttcore{A}{k + 1:d}}$, which is validated by the numerical results in \Cref{subsubsection:comparisonofhatt-1and-2}. This observation motivates our choice of using the singular values of $\VT{\ttcore{A}{1:k}}$ as weights. Unfortunately, we cannot establish a constant low bound for the singular values of $\VT{\ttcore{A}{1:k}}$, even though it arises from the  matricization of the Hadamard product of two left orthogonal tensors. This limitation prevents us from deriving a sharper result analogous to \Cref{lemma:lowerboundofA-B}. Nevertheless, the numerical results reported in \Cref{sec:experiments} show that the matrix $\HT{\ttcore{A}{k + 1:d}}$ is typically numerically low rank. Since $\mat{W}^{(k-1)}  = \HT{\ttcore{A}{k:d}}\left(\HT{\ttcore{R}{k:d}}\right)\tran$, this numerical low-rank property is inherited by $\mat{W}^{(k)}$. This argument is further supported by the numerical evidence presented in \Cref{sec:experiments}, although a rigorous theoretical proof remains challenging because $\HT{\ttcore{R}{k:d}}$ is not a standard random matrix. The potential numerical low-rank structure of $\mat{W}^{(k)}$ serves as the primary motivation for proposing the HPCRL-2 algorithm.

}

Furthermore, performing truncated SVD may reveal the singular values of $W^{(k)}$, which can help us select a more appropriate target rank $\ell_k$ to guarantee accuracy and efficiency.  
According to the computational complexity analysis provided in \autoref{sec:complexity-analysis}, HPCRL-2 is generally only slightly more complex than HPCRL-1. 
Therefore, HPCRL-2 may be beneficial in certain applications where highly accurate TT recompression is required.

\begin{figure}[!htpb]
  \centering
  \begin{tikzpicture}[baseline]
  \Vertex[x=-.5, y=0.75, size=0, Pseudo=True]{Ak-1}
  \Vertex[x=-.5, y=-0.75, size=0, Pseudo=True]{Rk-1}
  \Vertex[x=1, y=0.75, label={$ \ttcore{A}{k} $}, position={above}, distance=-1.2mm]{Ak}
  \Vertex[x=1, y=-0.75, label={$ \ttcore{R}{k} $}, position={below}]{Rk}
  \Vertex[x=2, y=0, label={$ \mat{W}^{(k)} $}, position={right}, distance=-1.2mm]{Wk}
  \Edge(Ak-1)(Ak)
  \Text[x=-.1, y=.75, fontsize=\tiny, color=black!75, position={below}, distance=2pt]{$ r_{k - 1}s_{k - 1} $}
  \Edge[label=$ \ell_{k - 1} $, position={above}, distance=.4](Rk-1)(Rk)
  \Edge[label=$ r_{k}s_{k} $, position={above right}](Ak)(Wk)
  \Edge[label=$ \ell_{k} $, position={below right}](Rk)(1.46666666666, -.4)
  \Edge[label=$ \ell_{k} $, position={below right}](Wk)(1.53333333333, -.35)
  \Edge[label=$ n_{k} $, position={left}](Ak)(1, .1)
  \Edge[label=$ n_{k} $, position={left}, distance=.4](Rk)(1, -.1)
  \draw[rounded corners, dashed] (1, 1.5) -- ++(-.6, -.8) -- ++(1.8, -1.35) -- ++(.6, .8) -- cycle;
  \draw[rounded corners, dashed] (.5, -.3) rectangle (1.5, -1.25);
  \draw[->] (2.2, .75) -- node[above] {\tiny Property \eqref{eq:kronecker_product_times_vec}} (3, .75);
  \draw[->] (1.7, -.75) -- node[below] {\parbox{1.3cm}{\tiny Arrange $ \ttcore{R}{k}(i_{k}) $ horizontally}} (3, -.75);
\end{tikzpicture} 
\begin{tikzpicture}[baseline]
  \Vertex[x=1, y=0.75, label={$ \mat{W}_{L} $}, position={above}]{WL}
  \Vertex[x=1, y=-0.75, label={$ \mat{W}_{R} $}, position={below}]{WR}
  \Edge(0, 0.75)(WL)
  \Edge[label = {$ n_{k}\ell_{k} $ }, position = {right}](1, .1)(WL)
  \Text[x=.1, y=.75, fontsize=\tiny, color=black!75, position=below, distance=2pt]{$ r_{k - 1}s_{k - 1} $}
  \Edge[label=$ \ell_{k - 1} $, position={above}](0, -0.75)(WR)
  \Edge[label = {$ n_{k}\ell_{k} $ }, position = {right}](1, -.1)(WR)
\end{tikzpicture} $ \longrightarrow $ 
\begin{tikzpicture}[baseline]
  \Vertex[x=1, y=0.75, label={$ \mat{W}_{L} $}, position={above}]{WL}
  \Vertex[x=1, y=-0.75, label={$ \mat{W}_{R} $}, position={below}]{WR}
  \Edge(0, 0.75)(WL)
  \Text[x=.1, y=.75, fontsize=\tiny, color=black!75, position=below, distance=2pt]{$ r_{k - 1}s_{k - 1} $}
  \Edge[label=$ \ell_{k - 1} $, position={above}](0, -0.75)(WR)
  \Edge[label = {$ n_{k}\ell_{k} $ }, position = {right}, color={red}](WL)(WR)
\end{tikzpicture}
  \caption{The process of HPCRL for direct representation \cref{HPCRL--1}.}\label{fig:HPCRL-No-SVD}
\end{figure}

\begin{figure}[!htpb]
  \centering
  \begin{tikzpicture}[baseline]
  \Vertex[x=-.5, y=0.75, size=0, Pseudo=True]{Ak-1}
  \Vertex[x=-.5, y=-0.75, size=0, Pseudo=True]{Rk-1}
  \Vertex[x=1, y=0.75, label={$ \ttcore{A}{k} $}, position={above}]{Ak}
  \Vertex[x=1, y=-0.75, label={$ \ttcore{R}{k} $}, position={below}]{Rk}
  \Vertex[x=1.75, y=0, label={$ \mat{W}^{(k)} $}, position={right}, distance=-1mm]{Wk}
  \Edge(Ak-1)(Ak)
  \node[below=2pt, color={black!75}] at (-0.1, 0.75) {\tiny $ r_{k - 1}s_{k - 1} $};
  \Edge[label=$ \ell_{k - 1} $, position={above}](Rk-1)(Rk)
  \Edge[label=$ n_{k} $, position={left}](Ak)(1, .1)
  \Edge[label=$ n_{k} $, position={left}, distance=.4](Rk)(1, -.1)
  \Edge[label=$ r_{k}s_{k} $, position={above right}](Ak)(Wk)
  \Edge[label=$ \ell_{k} $, position={below right}](Rk)(Wk)
\end{tikzpicture} \!\!$ \longrightarrow $ \!\!
\begin{tikzpicture}[baseline]
  \Vertex[x=1, y=0.75, label={$ \ttcore{A}{k} $}, position={above}, distance=-1mm]{Ak}
  \Vertex[x=1, y=-0.75, label={$ \ttcore{R}{k} $}, position=below, distance=-1mm]{Rk}
  \Vertex[x=2, y=0.75, label={$ \mat{U} $}, position=above]{U}
  \Vertex[x=2, y=-0.25, label={$ \mat{S} $}, position=right]{S}
  \Vertex[x=2, y=-0.75, label={$ \mat{V} $}, position=below]{V}
  \Edge(0, 0.75)(Ak)
  \Text[x=-.1, y=.75, fontsize=\tiny, color=black!75, position=below, distance=1.2pt]{$ r_{k - 1}s_{k - 1} $}
  \Edge[label=$ \ell_{k - 1} $, position={above}, distance=.2](0, -0.75)(Rk)
  \Edge[label=$ n_{k} $, position={left}](Ak)(1, .1)
  \Edge[label=$ n_{k} $, position={left}, distance=.4](Rk)(1, -.1)
  \Edge(Ak)(U)
  \Edge(U)(S)
  \Edge(S)(V)
  \Edge(V)(Rk)
  \draw[rounded corners, dashed] (0.5, 1.25) rectangle (2.5, 0.325);
  \draw[->] (2.8, 0.75) -- node[above] {\tiny Property \eqref{eq:kronecker_product_times_vec}} (3.6, 0.75);
  \draw[rounded corners, dashed, color = {red}] (0.5, -0.25) -- (1.5, -0.25) -- (1.5, 0) -- (2.5, 0) -- (2.5, -1.25) -- (0.5, -1.25) -- cycle;
  \draw[->] (2.8, -0.75) -- (3.6, -0.75);
\end{tikzpicture}
\begin{tikzpicture}[baseline]
  \Vertex[x=1, y=0.75, label={$ \mat{W}_{L} $}, position={above}]{WL}
  \Vertex[x=1, y=-0.75, label={$ \mat{W}_{R} $}, position={below}]{WR}
  \Edge(0, 0.75)(WL)
  \node[below=2pt, color={black!75}] at (0.1, 0.75) {\tiny $ r_{k - 1}s_{k - 1} $};
  \Edge[label=$ \ell_{k - 1} $, position={above}](0, -0.75)(WR)
  \Edge[label=$ n_{k}R $, position={right}, color = {red}](WL)(WR)
\end{tikzpicture}
  \caption{The process of HPCRL for (truncated) SVD \cref{HPCRL--2}.}\label{fig:HPCRL-process}
\end{figure} 

\begin{algorithm}[!htpb]
  \caption{HPCRL: Generate partial contractions $ \{\mat{W}^{(k)}\}_{k = 1}^{d - 1} $ of the Hadamard product $\tensor{Y} \odot \tensor{Z}$}\label{alg:HPCRL-main}
  \begin{algorithmic}[1]
    \Require Two TT tensors $ \tensor{Y} $ and $ \tensor{Z} $. A random TT tensor $ \tensor{R} $.
    \Ensure A set of sketching matrices $ \{\mat{W}^{(k)}\}_{k = 1}^{d - 1} $.
    \Function{$ \{\mat{W}^{(k)}\}_{k = 1}^{d - 1} $}{HPCRL}{$ \tensor{Y} $, $ \tensor{Z} $, $ \tensor{R} $}
    \State $ \mat{W}^{(d - 1)} = \HT{\ttcore{Y}{d} \pkp^{1} \ttcore{Z}{d}}\left(\HT{\ttcore{R}{d}}\right)\tran $ 
    \For{$ k = d - 1 $ down to $ 2 $}
      
      \State $ \mat{W}^{(k)} :=\mat{U}^{(k)}\mat{S}^{(k)}\left(\mat{V}^{(k)}\right)\tran$ 
      \For{$ i_{k} = 1 $ to $ n_{k} $}
        \For{$ \gamma = 1 $ to the number of non-zero elements in $\mat{S}^{(k)}$ }
          \State $ \mat{U}^{(k)}_{\gamma} = \mat{U}^{(k)}(:, \gamma)\big|^{s_{k}}_{r_{k}} $ 
          \State $ \mat{W}_{L}(:, \overline{i_{k}\gamma}) = \alg{vec}\big(\ttcore{Z}{k}(i_{k}) \mat{U}^{(k)}_{\gamma} (\ttcore{Y}{k}(i_{k}))\tran\big) $ \label{algline:WL} \Comment{$ \R^{s \times s} \times \R^{s \times r} \times \R^{r \times r} $}
          \State $ \mat{W}_{R}(:, \overline{i_{k}\gamma}) = \ttcore{R}{k}(i_{k})\mat{V}^{(k)}(:, \gamma) $ \label{algline:WR} \Comment{$ \R^{\ell \times \ell} \times \R^{\ell} $}
        \EndFor
      \EndFor
      \State $ \mat{W}^{(k-1)} = \mat{W}_{L}(\mat{I}_{n_{k}} \otimes \mat{S}^{(k)}) \mat{W}_{R}\tran $ 
    \EndFor
    \EndFunction
  \end{algorithmic}
\end{algorithm}

\begin{figure}[!htpb]
  \centering
  \resizebox{\textwidth}{!}{
  \begin{tikzpicture}
  \Vertex[x=0, y=0, label=$ \ttcore{A}{1} $, position=above, distance=-1mm]{A1}
  \Vertex[x=1, y=0, label=$ \ttcore{A}{2} $, position=above, distance=-1mm]{A2}
  \Vertex[x=2, y=0, label=$ \ttcore{A}{3} $, position=above, distance=-1mm]{A3}
  \Vertex[x=3, y=0, label=$ \ttcore{A}{4} $, position=above, distance=-1mm]{A4}
  \Vertex[x=0, y=-1.5, label=$ \ttcore{R}{1} $, position=below, distance=-1mm]{R1}
  \Vertex[x=1, y=-1.5, label=$ \ttcore{R}{2} $, position=below, distance=-1mm]{R2}
  \Vertex[x=2, y=-1.5, label=$ \ttcore{R}{3} $, position=below, distance=-1mm]{R3}
  \Vertex[x=3, y=-1.5, label=$ \ttcore{R}{4} $, position=below, distance=-1mm]{R4}
  \Edge[label=$ n_{1} $, position=right](A1)(R1)
  \Edge[label=$ n_{2} $, position=right](A2)(R2)
  \Edge[label=$ n_{3} $, position=right](A3)(R3)
  \Edge[label=$ n_{4} $, position=right](A4)(R4)
  \Edge[label=$ r_{1}s_{1} $, position=below](A1)(A2)
  \Edge[label=$ r_{2}s_{2} $, position=below](A2)(A3)
  \Edge[label=$ r_{3}s_{3} $, position=below](A3)(A4)
  \Edge[label=$ \ell_{1} $, position=above](R1)(R2)
  \Edge[label=$ \ell_{2} $, position=above](R2)(R3)
  \Edge[label=$ \ell_{3} $, position=above](R3)(R4)
  \draw[rounded corners, dashed] (2.5, 0.5) rectangle (4.3, -2) node[above left = -2pt] {\small $ \mat{W}^{(3)} $};
\end{tikzpicture}
\begin{tikzpicture}
  \Vertex[x=0, y=0, label=$ \ttcore{A}{1} $, position=above, distance=-1mm]{A1}
  \Vertex[x=1, y=0, label=$ \ttcore{A}{2} $, position=above, distance=-1mm]{A2}
  \Vertex[x=2, y=0, label=$ \ttcore{A}{3} $, position=above, distance=-1mm]{A3}
  \Vertex[x=0, y=-1.5, label=$ \ttcore{R}{1} $, position=below, distance=-1mm]{R1}
  \Vertex[x=1, y=-1.5, label=$ \ttcore{R}{2} $, position=below, distance=-1mm]{R2}
  \Vertex[x=2, y=-1.5, label=$ \ttcore{R}{3} $, position=below, distance=-1mm]{R3}
  \Vertex[x=3, y=-.75, label=$ \mat{W}^{(3)} $, position=right, distance=-1mm]{W3}
  \Edge[label=$ n_{1} $, position=right](A1)(R1)
  \Edge[label=$ n_{2} $, position=right](A2)(R2)
  \Edge[label=$ n_{3} $, position=right](A3)(R3)
  \Edge[label=$ r_{1}s_{1} $, position=below](A1)(A2)
  \Edge[label=$ r_{2}s_{2} $, position=below](A2)(A3)
  \Edge[label=$ \ell_{1} $, position=above](R1)(R2)
  \Edge[label=$ \ell_{2} $, position=above](R2)(R3)
  \Edge[label=$ r_{3}s_{3} $, position=above right](A3)(W3)
  \Edge[label=$ \ell_{3} $, position=below right](R3)(W3)
  \draw[rounded corners, dashed] (1.5, 0.5) rectangle (4, -2) node[above left = -2pt] {\small $ \mat{W}^{(2)} $};
\end{tikzpicture}
\begin{tikzpicture}
  \Vertex[x=0, y=0, label=$ \ttcore{A}{1} $, position=above, distance=-1mm]{A1}
  \Vertex[x=1, y=0, label=$ \ttcore{A}{2} $, position=above, distance=-1mm]{A2}
  \Vertex[x=0, y=-1.5, label=$ \ttcore{R}{1} $, position=below, distance=-1mm]{R1}
  \Vertex[x=1, y=-1.5, label=$ \ttcore{R}{2} $, position=below, distance=-1mm]{R2}
  \Vertex[x=2, y=-.75, label=$ \mat{W}^{(2)} $, position=right, distance=-1mm]{W2}
  \Edge[label=$ n_{1} $, position=right](A1)(R1)
  \Edge[label=$ n_{2} $, position=right](A2)(R2)
  \Edge[label=$ r_{1}s_{1} $, position=below](A1)(A2)
  \Edge[label=$ \ell_{1} $, position=above](R1)(R2)
  \Edge[label=$ r_{2}s_{2} $, position=above right](A2)(W2)
  \Edge[label=$ \ell_{2} $, position=below right](R2)(W2)
  \draw[rounded corners, dashed] (0.5, 0.5) rectangle (3, -2) node[above left = -2pt] {\small $ \mat{W}^{(1)} $};
\end{tikzpicture}}
  \caption{Generate sketches $ \{\mat{W}^{(k)}\}_{k = 1}^{3} $ for $ d = 4 $.}\label{fig:randorth-generate-W}
\end{figure}



\subsection{HaTT}
\label{QRforeachcore}
To develop an efficient TT recompression algorithm for the Hadamard product,  
we first generate a random Gaussian TT tensor $ \tensor{R}$ with the target TT ranks $ \{\ell_{k}\}_{k = 0}^{d} $, as defined in \cref{def:randTT-2}. Next, we compute partial contraction of TT tensors $\tensor{A}$ and $\tensor{R}$ to obtain the sketches $\mat{W}^{(k)}$ of $\tensor{A}$. These sketches $\mat{W}^{(k)}$ are sequentially defined in \cref{rtlpc} and can be efficiently computed using the HPCRL algorithm proposed in the previous subsection.  The process of generating sketches $\mat{W}^{(k)} $  is shown in \cref{fig:randorth-generate-W}.
Finally, inspired by the RandOrth method, starting with $k=1$ and initializing $\ttcore{X}{1} = \ttcore{A}{1}$, we sequentially construct a left-orthogonal compressed TT tensor $\tensor{X}$ as follows:
\begin{enumerate}
\item[(1)] Compute the QR factorization of the sketched matrix, i.e., 
\begin{equation}
\label{eq:updatcores1}
    \mat{X}^{(k)}{\mat{W}^{(k)}}=\mat{Q}^{(k)}\mat{R}^{(k)}. 
\end{equation}
\item[(2)] Update the cores of $\tensor{X}$ by 
\begin{equation}
\label{eq:updatcores2}
    \HT{\ttcore{X}{k+1}}=\mat{M}^{(k)}\HT{\ttcore{A}{k+1}} \quad \textnormal{and}\quad \VT{\ttcore{X}{k}}=\mat{Q}^{(k)}.
\end{equation}
\end{enumerate}

Note that $\mat{X}^{(k)}:=\VT{\ttcore{X}{k}}\in\R^{\ell_{k-1}n_k\times r_ks_k}$, $\mat{Q}^{(k)}\in\R^{\ell_{k-1}n_k\times \ell_k}$, $\mat{M}^{(k)}:=(\mat{Q}^{(k)})\tran\mat{X}^{(k)}\in\R^{\ell_k\times r_{k}s_{k}}$ and $\HT{\ttcore{A}{k+1}}\in\R^{r_ks_k\times n_{k+1}r_{k+1}s_{k+1}}$. In \cref{eq:updatcores2}, a  matrix-matrix multiplication of dimensions $ \R^{\ell_{k} \times r_{k}s_{k}} \times \R^{r_{k}s_{k} \times n_{k + 1}r_{k + 1}s_{k + 1}} $  must be performed. This operation incurs a computational cost of $(2rs - 1)nrs\ell$ floating-point operations and requires a storage space of $nr^2s^2$.  
In order to reduce the computational and memory complexity, we first rewrite $ \HT{\ttcore{X}{k + 1}} = \mat{M}^{(k)}\HT{\ttcore{A}{k + 1}} $ in TT tensor form as $\ttcore{X}{k + 1} = \mat{M}^{(k)} \times^{1} \ttcore{A}{k + 1}$.
The $i$th slice of TT core $\ttcore{X}{k + 1}$ is defined as $\ttcore{X}{k + 1}(i) = \mat{M}^{(k)}\ttcore{A}{k + 1}(i) = \mat{M}^{(k)}\big(\ttcore{Y}{k + 1}(i) \otimes \ttcore{Z}{k + 1}(i)\big)$.
According to property \cref{eq:kronecker_product_times_vec}, the $\gamma$th row of $\ttcore{X}{k + 1}(i)$ can be efficiently computed by:
\begin{equation}    
\begin{aligned}
\label{eq:M*core0}
  \ttcore{X}{k + 1}(\gamma, i, :) & = \mat{M}^{(k)}(\gamma, :)\big(\ttcore{Y}{k + 1}(i) \otimes \ttcore{Z}{k + 1}(i)\big) \\
  & = \left[\alg{vec}\bigl( \left(\ttcore{Z}{k + 1}(i)\right)\tran \mat{M}^{(k)}(\gamma, :)\big|^{s_{k}}_{r_{k}} \ttcore{\tensor{\tensor{Y}}}{k + 1}(i) \bigr)\right]\tran.
\end{aligned}
\end{equation}
Using \cref{eq:M*core0}, we can compute $ \ttcore{X}{k + 1} $ with a computational cost of $2nrs\ell(r + s - 1)$ flops. This approach eliminates the need for explicit computation and storage of the Hadamard product TT tensor core $\ttcore{\tensor{A}}{k}$, thereby reducing the overall complexity by an order of magnitude. In conclusion, the process of HaTT for recompressing $\tensor{Y} \odot \tensor{Z}$ is summarized in \cref{alg:HaTT}. Since there are two versions of HPCRL, we also obtain two versions of HaTT, denoted as HaTT-1 and HaTT-2, respectively.

\begin{algorithm}[!h]
  \caption{HaTT: Recompress $ \tensor{Y} \odot \tensor{Z} $ avoiding explicit representation of Hadamard product.}\label{alg:HaTT}
  \begin{algorithmic}[1]
    \Require Two TT tensors $ \tensor{Y}, \tensor{Z} $ with ranks $ \{r_{k}\}_{k = 0}^{d} $ and $ \{s_{k}\}_{k = 0}^{d} $. Target TT ranks $ \{\ell_{k}\}_{k = 0}^{d} $. 
    \Ensure A TT tensor $ \tensor{X} $ with ranks $ \{\ell_{k}\}_{k = 0}^{d} $
    \Function{$ \tensor{X} $}{\alg{HaTT}}{$ \tensor{Y} $ , $ \tensor{Z} $, $ \{\ell_{k}\}_{k = 0}^{d} $}
      \State Select a random Gaussian TT tensor $ \tensor{R} $ with target TT ranks $ \{\ell_{k}\}_{k = 0}^{d} $.
      \LineComment{Generate the sketching matrices $ \{\mat{W}^{(k)}\}_{k = 1}^{d - 1} $ from right to left}
      \State $ \{\mat{W}^{(k)}\}_{k = 1}^{d - 1} = \alg{HPCRL}(\tensor{Y}, \tensor{Z}, \tensor{R}) $ 
      \State $ \ttcore{X}{1} = \ttcore{Y}{1} \pkp^{3} \ttcore{Z}{1} $ 
      \For{$ k = 1 $ to $ d - 1 $}
        \State $ \mat{X}^{(k)} = \VT{\ttcore{X}{k}} $
        \State $ [\VT{\ttcore{X}{k}}, \sim] = \alg{QR}(\mat{X}^{(k)}\mat{W}^{(k)}) $
        \State $ \mat{M}^{(k)} = (\VT{\ttcore{X}{k}})\tran \mat{X}^{(k)} $  
        \LineComment{Contract $ \mat{M}^{(k)} $ onto $ \ttcore{Y}{k + 1} \pkp^{1, 3} \ttcore{Z}{k + 1} $ to get new $ \ttcore{X}{k + 1} $}
        \For{$ i = 1 $ to $ n_{k + 1} $}
          \For{$ \gamma = 1 $ to the number of rows of $ \mat{M}^{(k)} $}
            \State $ \mat{M}^{(k)}_{\gamma} = \mat{M}^{(k)}(\gamma, :)\big|^{s_{k}}_{r_{k}} $
            \State $ \ttcore{X}{k + 1}(\gamma, i, :) = \left[\alg{vec}\Big(\left(\ttcore{Z}{k + 1}(i)\right)\tran\mat{M}^{(k)}_{\gamma}\ttcore{Y}{k + 1}(i)\Big)\right]\tran $
          \EndFor
        \EndFor
      \EndFor
    \EndFunction
  \end{algorithmic}
\end{algorithm}

\section{Complexity analysis}
\label{sec:complexity-analysis}


Let us introduce the computational complexity of relevant operations as follows. For a matrix $ \mat{X} \in \R^{m \times n} $ with $ m > n $, the economy-sized QR (econ-QR) decomposition of $ \mat{X} $ is given by $ \mat{X} = \mat{Q}\mat{R} $, where $ \mat{Q} \in \R^{m \times n} $ has orthonormal columns and $ \mat{R} \in \R^{n \times n} $ is upper triangular. The econ-QR decomposition can be efficiently implemented using the Householder algorithm \cite{golub2013matrix}, which requires $4mn^2 - \frac{4n^3}{3}$ flops if only $\mat{Q}$ is needed.
The truncated SVD of $ \mat{X} $ is expressed as $ \mat{X} = \mat{U}\mat{S}\mat{V}\tran $, where $ \mat{U} \in \R^{m \times n} $ and $ \mat{V} \in \R^{n \times n} $ have orthonormal columns, and $ \mat{S} \in \R^{n \times n} $ is diagonal matrix. According to \cite{golub2013matrix}, computing the truncated SVD costs $O(mn^2)$ flops. 
Furthermore, the computational cost of matrix-matrix multiplication between $\mat{X} \in \R^{m \times n}$ and $\mat{Y} \in \R^{n \times r}$ is given by $m(2n - 1)r$ flops.
For convenience, we assume that the TT ranks of $ \tensor{Y} $, $ \tensor{Z} $, and $ \tensor{R} $ are $\{1,r,r,\ldots,r,1\}$, $\{1,s,s,\ldots,s,1\}$, and $\{1,\ell,\ell,\ldots,\ell,1\}$, respectively. Under these assumptions, the sketching matrices $ \{\mat{W}^{(k)}\}_{k = 1}^{d - 1} $ have dimensions $ \R^{rs \times \ell} $ and $ \varepsilon $-rank of $ R $. 



\cref{alg:HaTT} comprises two main components: the HPCRL algorithm for generating sketched matrices and the QR factorization of these sketched matrices for each core. In the remainder of this section, we analyze the computational complexity and memory requirements of both \alg{HPCRL} and \alg{HaTT}. 
In \alg{HPCRL}, we first construct the matrices $ \mat{W}_{L} $ and $ \mat{W}_{R} $, followed by computing $ \mat{W}_{L}(\mat{I} \otimes \mat{S}^{(k)}) \mat{W}_{R}$ to generate the sketched matrix $\mat{W}^{(k-1)}$. The complexity analysis for the two versions of HPCRL is presented below.

\emph{Complexity of HPCRL-1}. 
In HPCRL-1, the decomposition defined by \eqref{HPCRL--1} incurs no additional computational cost. 
To construct the matrix $ \mat{W}_{L}$, we perform $n_k\ell_k$ times of matrix products $\ttcore{Z}{k}(i_{k}) \mat{U}^{(k)}_{\gamma} (\ttcore{Y}{k}(i_{k}))\tran$. Here $\gamma=1,\,\ldots,\,\ell_k$, and $i_k=1,\,\ldots,\,n_k$, respectively. Given that $\ttcore{Z}{k}(i_{k})\in \mathbb{R}^{s\times s}$, $\mat{U}^{(k)}_{\gamma}\in \mathbb{R}^{s\times r}$, and $(\ttcore{Y}{k}(i_{k}))\tran\in \mathbb{R}^{r\times r}$, the computational cost for these matrix products is $n\cdot rs\ell(2r + 2s - 2)$. 
{\color{black} Constructing $ \mat{W}_{R} $ only costs some reshape, without matrix multiplication. Finally, multiplying $ \mat{W}_{L}\mat{W}_{R}\tran $ costs $ 2nrs\ell^2 - rs\ell $ flops.} This entire process is repeated for $d - 2$ iterations. Therefore, the total computational cost of \alg{HPCRL-1} is $(d - 2) \cdot [nrs\ell(2r + 2s + 2\ell)]$ flops.


\emph{Complexity of HPCRL-2}. 
Given that the matrix $ \mat{W}^{(k)} $ is of size $ {r_ks_k \times\ell_k} $, performing (truncated) SVD on it costs $ O(rs\ell^{2}) $. Assuming that we select the $ R $ largest singular values, $ \mat{W}_{L}\in \mathbb{R}^{r_{k-1}s_{k-1}\times n_kR}$ is constructed by performing $ n_kR $ times of matrix products $\ttcore{Z}{k}(i_{k}) \mat{U}^{(k)}_{\gamma} (\ttcore{Y}{k}(i_{k}))\tran$, where $\gamma=1,\,\ldots,\,R$, and $i_k=1,\,\ldots,\,n_k$, respectively. Here, $\ttcore{Z}{k}(i_{k}) \in \mathbb{R}^{s \times s}$, $\mat{U}^{(k)}_{\gamma} \in \mathbb{R}^{s \times r}$, and $(\ttcore{Y}{k}(i_{k}))\tran \in \mathbb{R}^{r \times r}$. The computational cost for these matrix products is $n \cdot rsR(2r + 2s - 2)$ flops. 
Similarly, constructing the matrix $ \mat{W}_{R} \in \mathbb{R}^{\ell_{k-1}\times n_kR}$ involves a computational cost of $ n\cdot R\ell(2\ell - 1) $ flops. {\color{black} Finally, multiplying $ \mat{W}_{L}\mat{W}_{R}\tran $ requires $  2nrsR\ell - rs\ell $ flops}. 
This process is repeated $d - 2$ times. Therefore, the total computational cost of \alg{HPCRL-2} is approximately $ (d-2) \cdot [nrsR(2r+2s+2\ell)+O(rs\ell^2)] $.

\emph{Complexity of HaTT}. Let $ C_{\alg{HPCRL}}$  denote the computational cost of \alg{HPCRL}. After executing \alg{HPCRL}, the HaTT algorithm first multiplies the core  $ \mat{X}^{(k)} $ by $ \mat{W}^{(k)} $,  producing a matrix of size $ {\ell_{k-1}n_k \times \ell_{k}} $ matrix, which incurs a computational cost of $ n\ell^2(2rs-1) $ flops. 
Next, the econ-QR factorization of the resulting matrix of size $ {\ell_{k-1}n_k \times \ell_{k}} $ is performed, with a computational cost of $ (4n  - 4/3)\ell^{3} $. Following this, $ \mat{M}^{(k)} = (\mat{Q}^{(k)})\tran\mat{X}^{(k)} $ is computed, and it is used to update the core $ \ttcore{X}{k + 1} $ as \cref{eq:M*core0}. These steps cost $ rs\ell(2\ell n - 1) $ and $ nrs\ell(2r + 2s - 2) $ flops, respectively. 
By sweeping through the entire TT tensor once, the total computational cost of \alg{HaTT} is:
\begin{align}
\label{CHATT}
  C_{\alg{HaTT}} = C_{\alg{HPCRL}} + (d-2)[nrs\ell(2r + 2s + 4\ell)].
\end{align}

By substituting the computational cost of HPCRL into \eqref{CHATT}, we obtain the total computational cost of HaTT and summarize it, along with that of other TT recompression algorithms, in \cref{tab:TT-Rounding-flops}.
{\color{black}
Compared with other TT recompression algorithms, the first advantage of the HaTT algorithm lies in the zero-cost evaluation of the Hadamard product, which effectively avoids both computational and storage bottlenecks. For $\ell \leq rs$, the advantage of HaTT-1 in terms of recompression cost over other TT recompression algorithms is summarized as follows. 
\begin{itemize}
    \item Compared with TT-Rounding: The speedup of HaTT-1 is given by $$\text{SP}(\ell):=\frac{5r^2s^2+6rs\ell+2\ell^2}{4r\ell+4s\ell+6\ell^2},$$
    which is a strictly decreasing function with respect to $\ell\in[0,rs]$. Therefore, the minimal speedup of HaTT-1, achieved at $\ell=rs$, is $\frac{13}{6+4/s+4/r}$, which is greater than 1 for $r+s\geq 4$. 
    In the special case $\ell=r=s$, the speedup of HaTT-1 simplifies to $\frac{5r^2+6r+2}{14}$, which is of the order of $r^2$.
    \item Compared with OrthRand: The speedup of HaTT-1 is given by $$\text{SP}(\ell):=\frac{5r^2s^2+2rs\ell+4\ell^2}{4r\ell+4s\ell+6\ell^2},$$
    which exhibits a similar behavior to that in the TT-Rounding case. Therefore, the detailed analysis is omitted here.
    \item Compared with RandOrth: The speedup of HaTT-1 is given by $$\text{SP}(\ell):=\frac{4rs+6\ell}{4r+4s+6\ell},$$
    which is a strictly decreasing function with respect to $\ell\in[0,rs]$ and $\min\{r,\,s\}$ $\geq 2$. Therefore, the minimal speedup of HaTT-1, achieved at $\ell=rs$, is $\frac{10}{6+4/s+4/r}$, which is greater than 1 for $\min\{r,\,s\}$ $\geq 2$.
    In the special case $\ell=r=s$, the acceleration of HaTT-1 simplifies to $\frac{4r+6}{14}$, which scales linearly with $r$.
    The comparison between TwoSided and HaTT-1 can be carried out in a similar manner. 
\end{itemize}
In summary, HaTT-1 is consistently more efficient than other recompression algorithms. When the target rank $\ell$ is much smaller than the product of the original ranks $rs$ (i.e., $\ell \ll rs$), the advantage of HaTT is most pronounced. The speedup analysis for HaTT-2 can be carried out in a similar, albeit more complex, manner. Its speedup is very similar to that of HaTT-1, particularly for $\ell \ll rs$, since the leading term is the same.

}

\begin{table}[!htbp]
  \centering
  \caption{Summary of the computational costs (flops) of HaTT and other TT recompression algorithms. }\label{tab:TT-Rounding-flops}
  \begin{tabular}{ccl}\hline
    Algorithms & Hadamard product cost  & Recompression cost   \\\hline
    \alg{TT-Rounding} & $O(dnr^2s^2)$ & $ (d - 2)n(5r^{3}s^{3} + 6r^{2}s^{2}\ell + 2rs\ell^{2}) $ \\
    \alg{OrthRand} & $O(dnr^2s^2)$ & $ (d - 2)n(5r^{3}s^{3} + 2r^{2}s^{2}\ell + 4rs\ell^{2}) $ \\
    \alg{RandOrth} & $O(dnr^2s^2)$ & $ (d - 2)n(4r^{2}s^{2}\ell + 6rs\ell^{2}) $ \\
    \alg{TwoSided} & $O(dnr^2s^2)$ & $ (d - 2)n(6r^{2}s^{2}\ell + 6rs\ell^{2}) $ \\
    \hline
    \alg{HaTT-1} & 0 & $  (d-2)[nrs\ell(4r+4s+6\ell)]$
    \\\hline
   \multirow{2}{*} {\alg{HaTT-2}} & \multirow{2}{*} {0} & $  (d-2)[nrs(2r+2s+2\ell)(\ell+R)$
    \\
   & & $~~~~+2nrs\ell^2+O(rs\ell^2)]$\\\hline
  \end{tabular}
\end{table}

\section{Numerical experiments} \label{sec:experiments}

We conducted four experiments to evaluate the efficiency and accuracy of the \alg{HaTT} algorithm. The HaTT implementation includes two versions: HaTT-1, which directly represents $ \mat{W}^{(k)} $ using \cref{HPCRL--1}, and HaTT-2, which employs the SVD \cref{HPCRL--2} of $\bm{W}^{(k)}$. For comparison, we used four state-of-the-art algorithms as baselines: TT-Rounding \cite{oseledets2011tensor}, and three randomized algorithms (RandOrth, OrthRand, and TwoSided) proposed in \cite{al2023randomized}. 
All algorithms were implemented using the TT-Toolbox \cite{oseledets2011matlab}. Our numerical experiments were performed on a machine equipped with an AMD EPYC 7452 CPU and 256 GB of RAM, using MATLAB R2020a. Each simulation was run five times with different random seeds, and we report the mean value of errors and computational times. For efficiency comparison, the relative error was measured in all experiments to ensure that our method and the benchmark methods achieved comparable accuracy. The code for our experiments is available on GitHub: \url{https://github.com/syvshc/HaTT_code}.



\subsection{Example 1: Hadamard product of Fourier series functions}

To evaluate the performance of HaTT with different target ranks, we consider the Hadamard product of two Fourier series functions: 
$$y(t) = \sum_{j = 1}^{60} a_{j}\sin(jt),~~ \hbox{and}~~z(t) = \sum_{j = 1}^{60} b_{j}\cos(jt),$$
where $ \{a_{j},b_j:j=1,2,\ldots,60\} $ are random, independent, uniformly distributed values from  $[0.1, 10.1]$. The function values of $y(t)$ and $z(t)$ at $t_i=\frac{2\pi i}{20^7}$, with $i=1,2,\cdots,20^7$, are folded to two $7$th-order tensors $\tensor{Y}$ and $\tensor{Z}\in\mathbb{R}^{20\times20\times\cdots\times20}$, respectively. These tensors are then represented in the TT format using TT-SVD \cite{oseledets2011tensor}. We apply TT recompression algorithms to recompress the Hadamard product $\tensor{Y}\odot\tensor{Z}$. For all simulations, the target TT rank $\ell$ increases from $4$ to $60$ with the step size of $4$. {\color{black} For all target ranks, we set the truncation tolerance to \texttt{ep = 1e-5} in the HaTT-2 algorithm of $\mat{W}^{(k)}$.}

\begin{figure}[!hptb]
  \centering
  \includegraphics[width=0.95\textwidth]{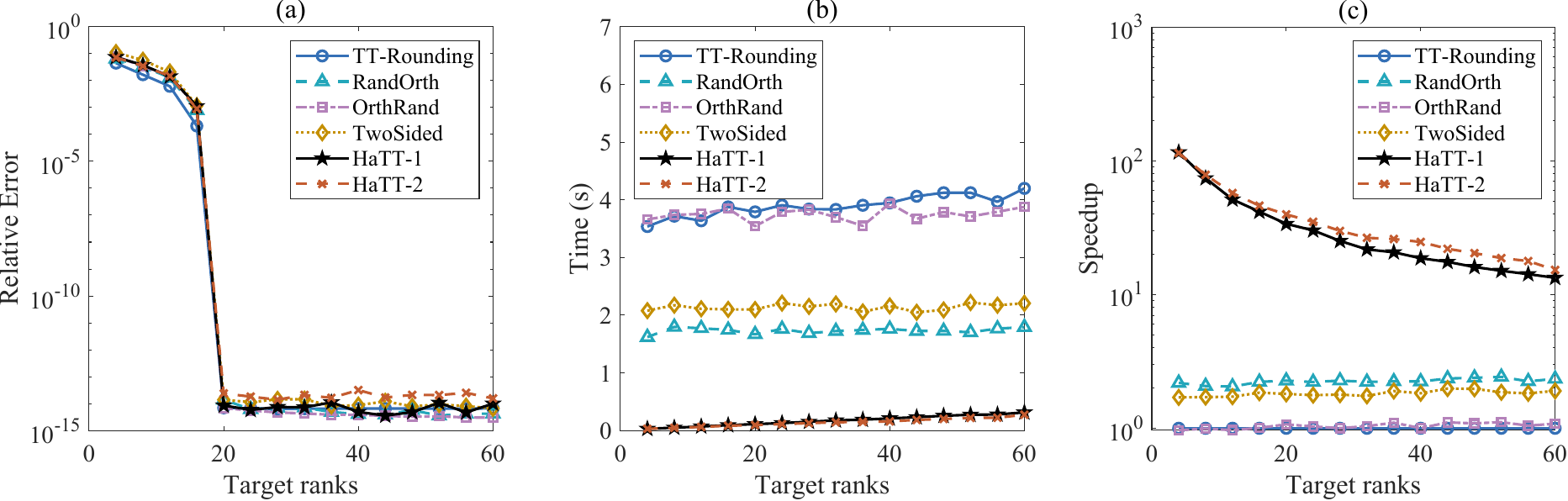}
  \caption{(a) Relative errors, (b) running times, and (c) speedups for Example 1.}\label{fig:sincos}
\end{figure}

We report the relative errors, running times, and speedups (compared to the TT-Rounding algorithm) for all TT recompression algorithms in \cref{fig:sincos}.
As shown in \cref{fig:sincos}(a), the relative errors obtained by \alg{HaTT} are comparable to those of other baseline methods, validating the accuracy of \alg{HaTT}.
From \cref{fig:sincos}(b), we observe that the randomized algorithms are faster than TT-Rounding, with \alg{HaTT} proving particularly effective due to the Hadamard-avoiding technique.
It follows from \cref{fig:sincos}(c) that the speedup of HaTT compared to other baseline methods decreases as the target TT rank $\ell$ increases.
For this example, the speedup of \alg{HaTT} compared to TT-Rounding, RandOrth, OrthRand, and TwoSided ranges {\color{black} from $13.3\times$ to $ 115.6\times$, $5.7\times$ to $52.9\times$, $12.3\times$ to $119.6\times$, and $7\times$ to $67.9\times$, respectively.} The relative error is defined as 
\begin{align*}
  \|\tensor{A}^{\mathit{round}} - \tensor{A}\|_{\mathrm{F}}/\fnorm{\tensor{A}},
\end{align*}
where $\tensor{A}^{\mathit{round}}$ is the low rank TT tensor obtained by TT recompression algorithm and $\tensor{A}$ is the exact solution of $\tensor{Y}\odot\tensor{Z}$. We can see that when target rank is larger than $15$, relative error of all methods is less than $10^{-13}$, which means that all methods, including \alg{HaTT}, can recompress the TT tensor with accuracy approaching machine precision.



\subsubsection{Comparison between HaTT-1 and HaTT-2}
\label{subsubsection:comparisonofhatt-1and-2}

{\color{black} We use Example 1 to provide a clear comparison between the performances of HaTT-1 and HaTT-2. When computing the truncated SVD of $ \mat{W}^{(k - 1)} $ in HaTT-2, the following truncation strategy is considered: retain only the singular values larger than $ \text{\texttt{ep}}\times \|\mat{W}^{(k - 1)}\|_{2} $. 
We mainly focus on how this truncation criterion affects the efficiency and accuracy of HaTT-2. The tensor to be recompressed, $\tensor{A} = \tensor{Y} \odot \tensor{Z}$, remains unchanged, and the target rank $\ell$ is set to 20. Tensors $\tensor{Y}$ and $\tensor{Z}$ are first recompressed using the TT-SVD algorithm \cite{oseledets2011tensor} with a relative tolerance of $\varepsilon_r = 10^{-14}$.
The resulting TT ranks of the two tensors $\tensor{Y}$ and $\tensor{Z}$, respectively, are $[1, 12, 37, 24, 14, 12, 20, 1] $ and $[1, 13, 35, 22, 14, 13, 20, 1]$. We then compute their Hadamard product $\tensor{A} = \tensor{Y} \odot \tensor{Z}$. As shown in \Cref{fig:SVD}, even the matrices $\unfold{A}{3}$ and $\unfold{A}{4}$ are not low rank, the $6$-th singular value of $\unfold{A}{3}$ and the $8$-th singular value of $\unfold{A}{4}$ both below $10^{-14} \norm{\unfold{A}{3}}_2$ and $ 10^{-14} \norm{\unfold{A}{4}}_2$, respectively. If we apply the TT-rounding algorithm to the Hadamard product $\tensor{A}$ with the same relative tolerance $\varepsilon_r = 10^{-14}$, yielding relative $10^{-14}$-numerical TT-ranks of $[1, 7, 7, 8, 11, 13, 20, 1] $, which is consistent with our theoretical analysis (see inequality \eqref{theorem3.5:eq-1}). These results indicate that the Hadamard product $\tensor{A}$ exhibits relatively low numerical ranks.

The singular value distributions of $\HT{\ttcore{A}{k + 1:d}}$ for $k = 3$ and $4$ are shown in \Cref{fig:SVD}, exhibiting a similarly rapid decay pattern to that of $\unfold{A}{k}$. It is clear that $\HT{\ttcore{A}{k + 1:d}}$ is also relative numerical low rank, though slightly higher than that of $ \unfold{A}{k}$. The singular value distributions of $\VT{\ttcore{A}{1:k}}$, also presented in \Cref{fig:SVD}, indicate that most singular values are on the order of $10^{-2}$. By choosing $j$ such that $\sigma_j\left(\VT{\ttcore{A}{1:k}}^+\right)\approx 10^{-2}$  (specifically, $j=40$), and according to \Cref{theorem3.5numericallowrankofH}, the relative ${10^{-16}}$-numerical rank of $\HT{\ttcore{A}{k + 1:7}}$ does not exceed $r_\varepsilon^k+j-1$. This observation is fully consistent with the theoretical prediction in \Cref{theorem3.5numericallowrankofH}, further confirming that $\HT{\ttcore{A}{k + 1:d}}$ is indeed numerically low rank.

Next, we use the HaTT-1 and HaTT-2 to recompress tensor $\tensor{A}$, with the low rank tesnors $\tensor{Y}$ and $\tensor{Z}$ obtained by TT-SVD. We set \texttt{ep} to 0.95, 1e-1, 4e-2, 1e-3, 1e-5, 1e-10, and 1e-14, respectively. For each trial under the same parameter setting, identical random seeds are used for both HaTT-1 and HaTT-2 to ensure a fair comparison. It is important to note that the results of HaTT-1 depend solely on the random seed, and are unaffected by the value of \texttt{ep}. The computation time and relative error of HaTT-1 and HaTT-2 are presented in \Cref{fig:timeerr}. As observed, HaTT-2 becomes faster than HaTT-1 as \texttt{ep} increases. Except for the case of $\texttt{ep}=0.95$, HaTT-2 produces satisfactory low-rank tensor approximations, achieving relative errors on the order of $10^{-10}$ or even smaller.

The singular value distributions of $W^{(3)}$  and $W^{(4)}$ with $\ell=20$ are also presented in \Cref{fig:SVD}, exhibiting a similarly trend as that observed for $\HT{\ttcore{A}{4:7}}$ and $\HT{\ttcore{A}{5:7}}$. It is evident that both $W^{(3)}$ and $W^{(4)}$ are numerically low rank. This serves as the primary motivation for proposing the HaTT-2 algorithm. Next, we investigate whether this numerically low rank property persists as $\ell$ varies. The singular value distributions of $W^{(3)}$  and $W^{(4)}$ with $\ell=40$ and $60$ are also presented in \Cref{fig:SVD}, clearly indicating that all these matrices remain numerically low rank. Moreover, we observe that the singular value distributions of the sketch matrices remain consistent as $\ell$ increases. These observations suggest that the singular value distribution provides valuable guidance for determining whether the target rank $\ell$ is sufficiently large to achieve an accurate low-rank TT approximation. For example, if the smallest singular value of sketch matrix is not as small as excepted, it indicates that $\ell$ should be increased. Otherwise, f the smallest singular value is significantly smaller than expected, $\ell$ may exceed the numerical rank of $\tensor{A}$, leading to unnecessary computational overhead. 

Since the target rank $\ell$ is typically unknown in practice, HaTT-1 must balance accuracy and efficiency by heuristically selecting $\ell$. In contrast, the HaTT-2 can achieve both high accuracy and efficiency even when $\ell$ is set to a relatively large value. To validate this claim, we fix \texttt{ep = 1e-3} and compare the performance of HaTT-1 and HaTT-2 with various target ranks $\ell=[10,\,20,\,30,\,40,\,50,\,60]$. The computation time and the relative error as functions of $\ell$ are shown in \Cref{fig:HaTT1and2——testrank}. Across all target ranks, HaTT-2 consistently outperforms HaTT-1 in computational efficiency while maintaining a comparable level of accuracy. Moreover, the computation time of HaTT-2 with $\ell=30$ is approximately the same as that of HaTT-1 with $\ell=20$, further supporting our claim.

\begin{figure}[h!]
    \centering
    \includegraphics[width=0.8\textwidth]{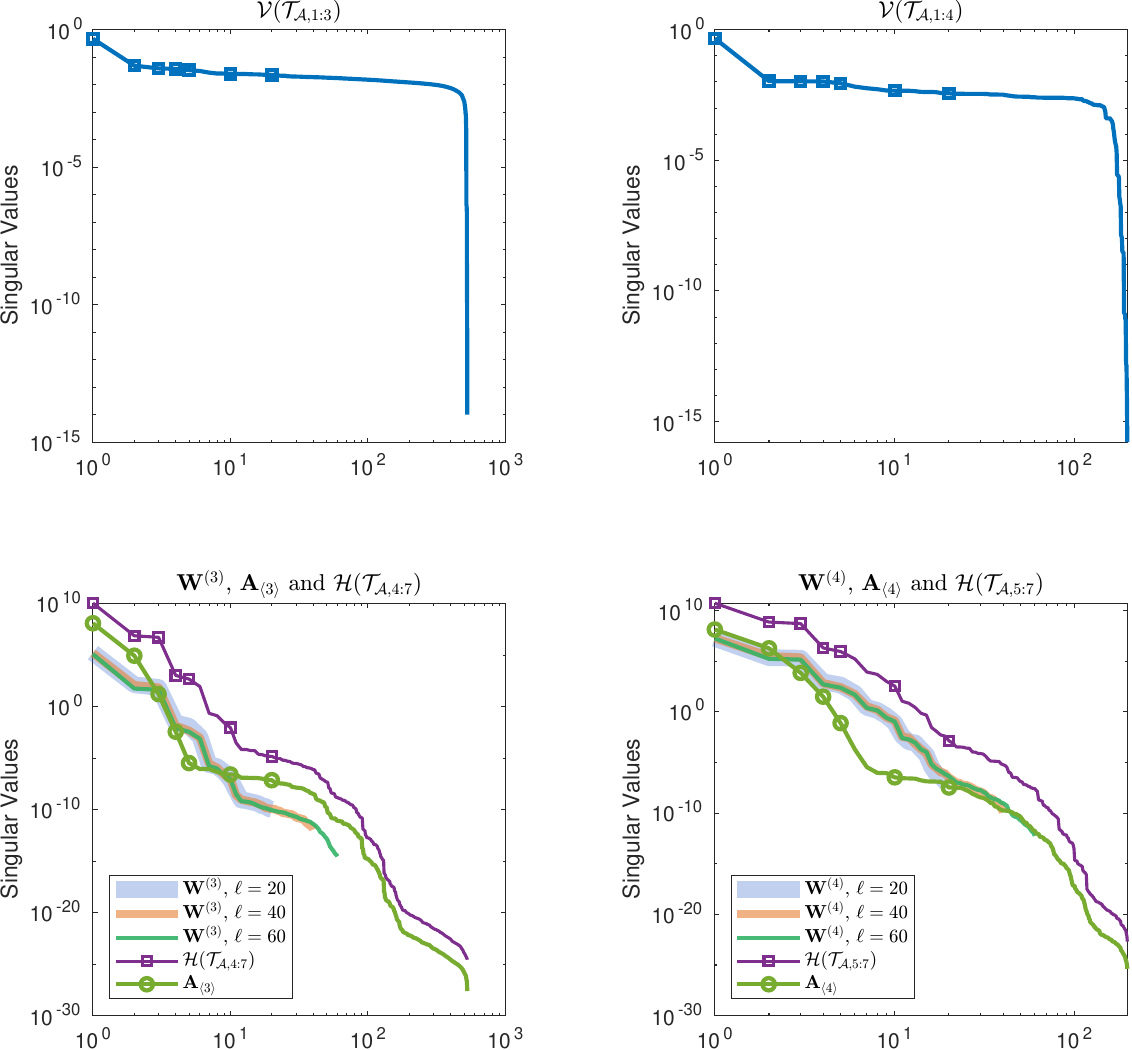}
    \caption{Singular values of $\VT{\ttcore{A}{1:k}}$, $\HT{\ttcore{A}{k + 1:d}}$, $ \mat{W}^{(k-1)} $, and $ \unfold{A}{k - 1} $ with $ k = 4, 5 $.}
    \label{fig:SVD}
\end{figure}

\begin{figure}[h!]
    \centering
    \includegraphics[width=0.9\textwidth]{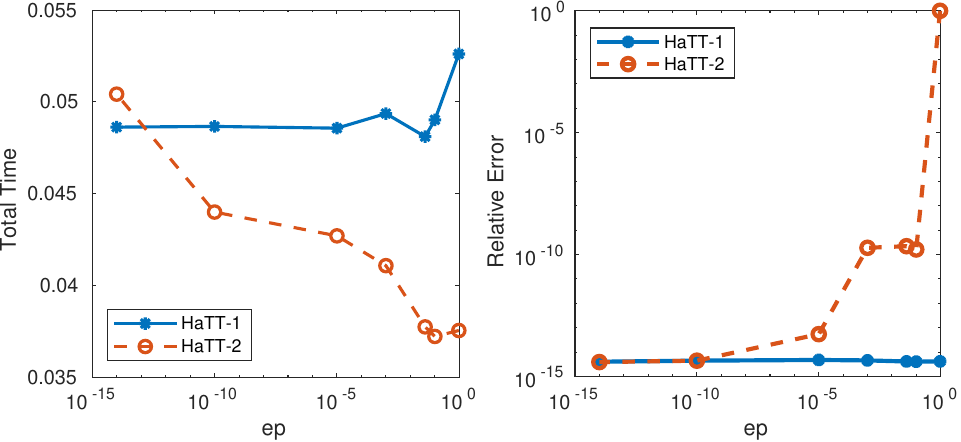}
    \caption{Comparison of HaTT-1 and HaTT-2 with varying values of \texttt{ep}. Here, $\ell = 20$. Left: computation time; Right: relative error.}
    \label{fig:timeerr}
\end{figure}

\begin{figure}
    \centering
    \includegraphics[width=0.9\linewidth]{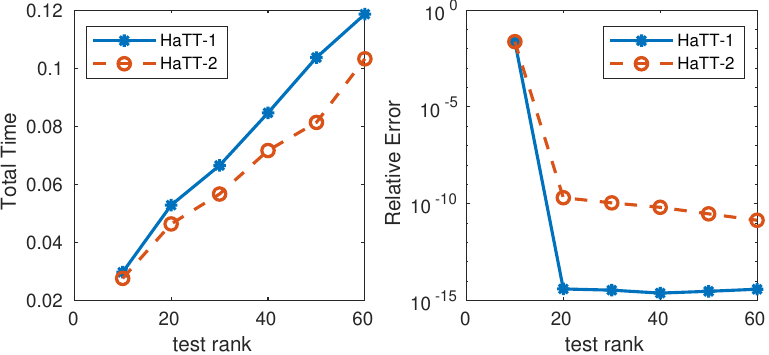}
    \caption{Total computation time (left) and relative error (right) for different target ranks $\ell $ with a fixed \texttt{ep = 1e-3}}
    \label{fig:HaTT1and2——testrank}
\end{figure}

}

\subsection{Example 2: Hadamard product of random TT tensors with different TT ranks}

In this example, we study the performance of HaTT for the Hadamard product of TT tensors with different ranks. Two $7$th-order TT tensors $\tensor{Y},\tensor{Z}\in\R^{20\times20\times\ldots\times20}$ are set as 
random uniform TT tensors (see \cref{def:randTT-2}). The corresponding maximal TT rank $r$ or $s$ increases from 60 to 150 in steps of 10. The target TT rank $\ell$ is fixed to $60$. The relative error in this experiment is defined in the same way as in the previous experiment.
The relative errors, running times, and speedups for all TT recompression algorithms are shown in \cref{fig:randtest_vary_r}. 
According to \cref{fig:randtest_vary_r}(a), the accuracy of HaTT is almost the same as RandOrth and OrthRand, though slightly larger than the accuracy of TT-Rounding.
As $r$ or $s$ increases, the speedup of HaTT compared to other baseline methods increases rapidly. For $r=s=110$, the HaTT
can achieve $6.9\times\sim 109.6\times$ speedup. 
It should be noted that simulations for $r$ or $s$ larger than $110$ using the baseline methods could not be performed due to memory limitations. Since HaTT avoids explicit representation of the Hadamard product $\tensor{Y}\odot\tensor{Z}$, it can recompress the Hadamard product even when $r=s=150$. This observation fully demonstrates the advantage of HaTT in terms of memory efficiency.

\begin{figure}[!htpb]
  \centering
  \includegraphics[width=0.95\textwidth]{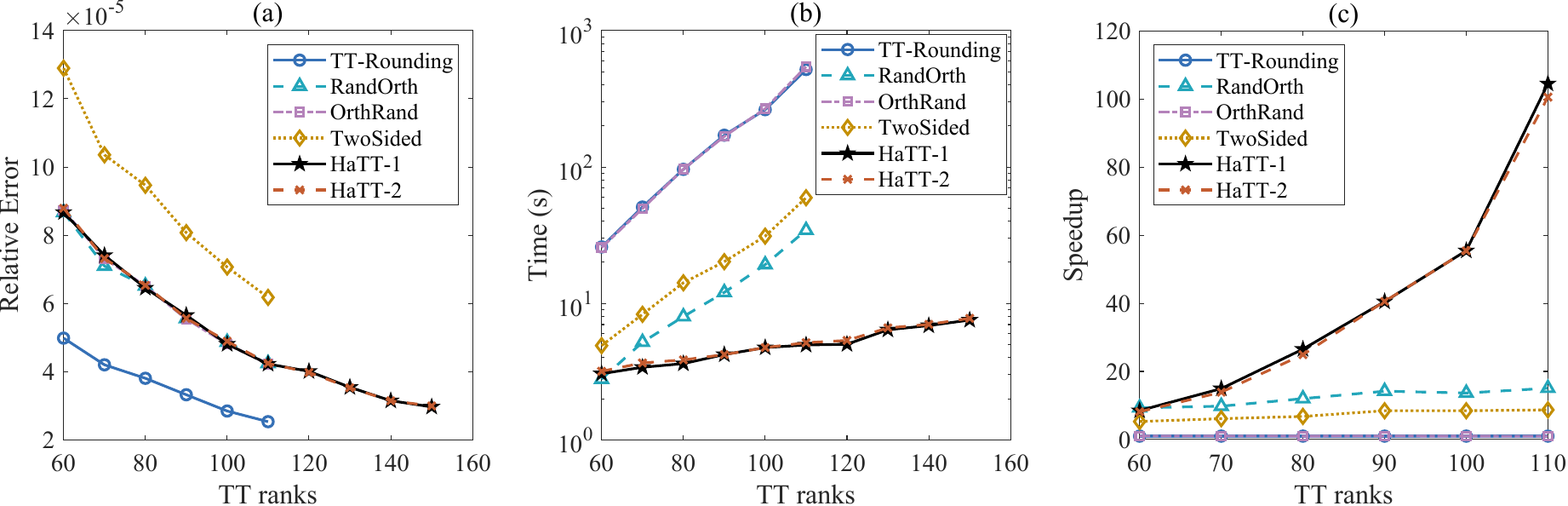}
  \caption{(a) Relative errors, (b) running times, and (c) speedups for Example 2.}\label{fig:randtest_vary_r}
\end{figure}



\subsection{Example 3: application of HaTT in power iteration}
We discuss the performance of HaTT in power iteration, which is applied to find the largest element of a TT tensor $\tensor{Y}$ \cite{espig2020iterative}. We generate the tensor $\tensor{Y}$ using a multivariate function, i.e., Qing or Alpine \cite[Func. 98 and 6]{jamil2013literature}
\begin{align*}
  \text{Qing function:}\ \ f(x_{1}, \dots, x_{d}) & = \sum_{i = 1}^{d} (x_{i} - i)^{2},\ x_i\in[-500,500],\\
  \text{Alpine function:}\ \ f(x_{1}, \dots, x_{d}) & = \sum_{i = 1}^{d} |x_{i}\sin(x_{i}) + 0.1x_{i}|,\ x_i\in[-2.5\pi,2.5\pi],
\end{align*}
where $d=10,\,20,\,\ldots$, or $50$ is the dimension. The function is discretized on a uniform grid with $10^d$ mesh points, resulting in tensor $\tensor{Y}\in\R^{10\times10\times\cdots\times 10}$. Since the functions are summations of $d$ separable multivariate functions, the tensor $\tensor{Y}$ can be represented as TT tensors with TT ranks
$ \{1, d, \dots, d, 1\} $. In the power iteration, the initial value is set to a tensor whose elements are all $1$, and the maximum number of iterations is set to $100$. We set the target TT rank $\ell$ to $5$ in all TT recompression algorithms, and \texttt{ep=1e-5} for HaTT-2. The relative errors of the largest elements obtained by the power iteration equipped with different recompression algorithms are displayed in \cref{fig:find_largest} (a) and (d), which indicates that the accuracy of the power iteration equipped with any recompression algorithm is acceptable. As shown in \cref{fig:find_largest} (b--c) and (e--f), the power iteration equipped with HaTT is much faster than other baseline methods, especially for larger $d$. 

\begin{figure}[!hptb]
  \centering
  \includegraphics[width=0.95\textwidth]{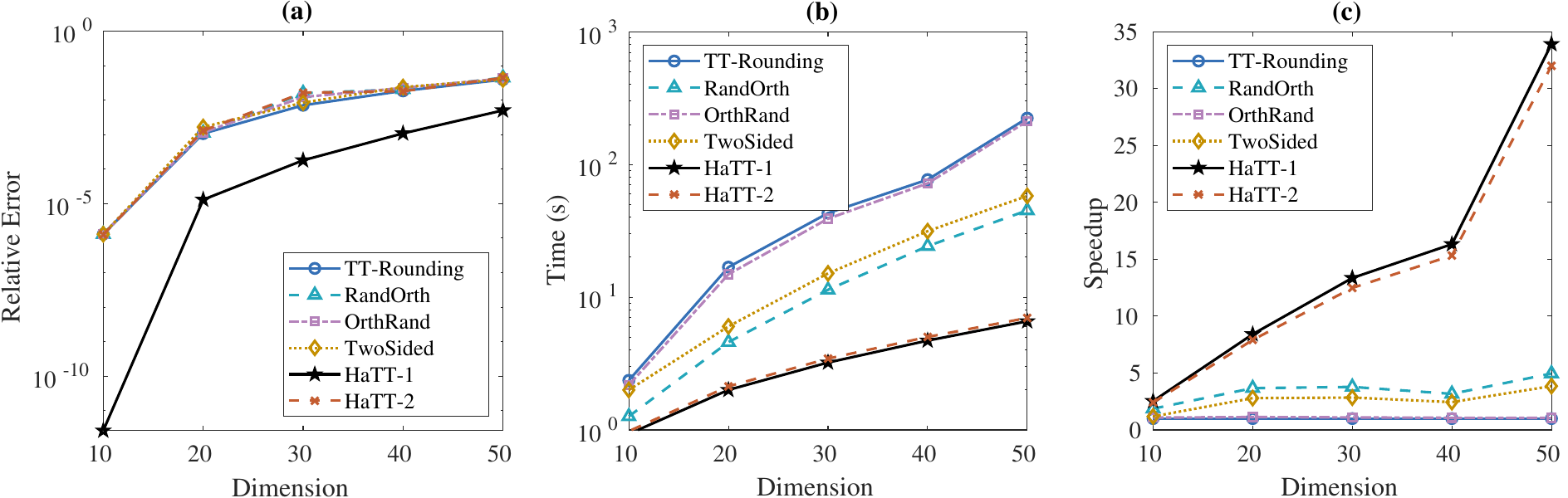} \\[1ex]
  \includegraphics[width=0.95\textwidth]{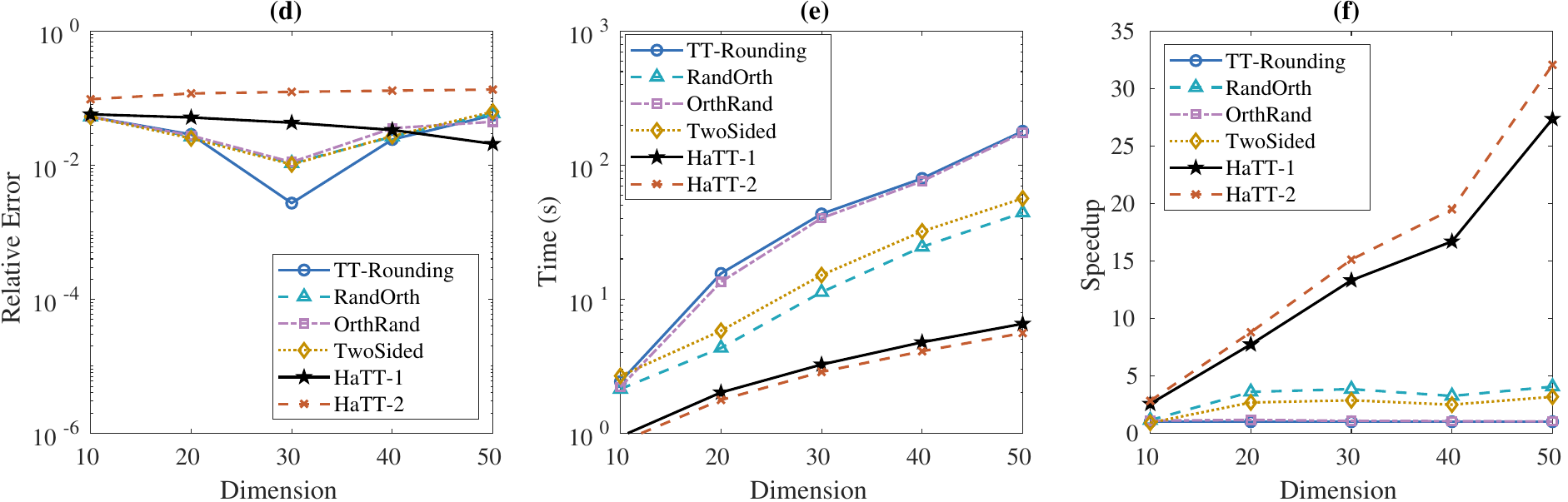}
  \caption{Example 3: power iteration equipped with different recompression algorithms for Qing (above) and Alpine (below) functions.}\label{fig:find_largest}
\end{figure}

\subsection{Example 4: application of HaTT in solving nonlinear PDEs}
We investigate the performance of HaTT in solving nonlinear PDEs. Let us consider the following Allen--Cahn equation with periodic boundary conditions:
\begin{equation}\left\{
\begin{aligned}
  & \frac{\partial\phi}{\partial t} = -\gamma\frac{\delta \mathcal{E}[\phi]}{\delta \phi}, \quad \hbox{in~} \Omega \times (0, T],\\
  &\phi = \phi_0,\quad \hbox{in~} \Omega,\label{eq:gradient-flow} 
\end{aligned}
\right.
\end{equation}
where $\phi:=\phi(\vb{x},t) $, $ \Omega \subset\mathbb{R}^2$ or $\mathbb{R}^3 $, $ \mathcal{E}[\phi] $ is a given total free energy functional, $\gamma$ is the mobility, and $ \frac{\delta \mathcal{E}[\phi]}{\delta \phi} $ is the variational derivative, respectively. The free energy functional usually can be written explicitly as: 
\begin{align}\label{eq:free-energy}
  \mathcal{E}[\phi] = \frac{1}{2} \int_{\Omega}| \nabla\phi|^2 \textnormal{d}\vb{x} + \int_{\Omega}\frac{1}{4\varepsilon^2}(1-\phi^2)^2 \textnormal{d}\vb{x},
\end{align}
where $\varepsilon$ is a small parameter. The Allen--Cahn equation \eqref{eq:gradient-flow}
satisfies the following energy dissipation law:
\begin{align*}
  \frac{\dd}{\dd t}\mathcal{E}[\phi] = \left\langle \frac{\delta \mathcal{E}}{\delta \phi}, \frac{\partial \phi}{\partial t} \right\rangle  =- \left\langle \frac{\delta \mathcal{E}}{\delta \phi}, \frac{\delta \mathcal{E}}{\delta \phi} \right\rangle \le 0,
\end{align*}
where $\left<\varphi,\,\psi\right>=\int_\Omega \varphi\psi\textnormal{d} \vb x$.
As is well known, several popular approaches exist for discretizing \eqref{eq:gradient-flow}, including the convex splitting method \cite{baskaran2013convergence, elliott1993global,  eyre1998unconditionally, shen2012second}, 
the stabilization method \cite{shen2010numerical, zhu1999coarsening}, the exponential time differencing \cite{ju2015fast, zhang2016extreme}, the invariant energy quadratization (IEQ) \cite{yang2016linear, zhao2017numerical}, the scalar auxiliary variable (SAV) \cite{shen2018scalar, shen2019new}, the discrete variational derivative (DVD) method \cite{du1991numerical, huang2020parallel, huang2022parallel,  PFCDVD, huang2024high}, all of which preserve the energy or modified energy dissipation law at the discrete level. In this work, we employ the relaxed DVD scheme \cite{huang2024high} for the time discretization and the second-order central difference scheme for the spatial discretization. In the relaxed DVD scheme \cite{huang2024high}, for $m\in\mathbb{Z}_+$, the free energy functional \eqref{eq:free-energy} is modified as:
\begin{align}\label{eq:modified-free-energy}
  \bar{\mathcal{E}}[\phi,\nu] =  \frac{1}{2}\int_{\Omega}\left(| \nabla\phi|^2+\frac{\beta}{\varepsilon^2}\phi^2 \right)\textnormal{d}\vb{x} + \nu^m,
\end{align}
where $\beta>0$ is a  suitable stabilization parameter and $\nu:=\nu(\phi)=\left(\int_\Omega\bar{E}(\phi) \textnormal{d}\vb{x}\right)^{1/m}$ with $\bar{E}(\phi)=\frac{1}{4\varepsilon^2}(1+\beta-\phi^2)^2$, respectively.

Assume that $[0,\,T]$ is divided into time intervals $[t^w,\,t^{w+1}]$, where $w=0,\,1,\,\ldots,$ is the time step index and $t^w=w\Delta t$ with $\Delta t$ being the time step size. Let us denote $\phi^w\approx\phi(\vb{x},t^w) $ and $\nu^w\approx\nu(\vb{x},t^w) $. The discrete variational derivative in the relaxed DVD scheme is defined as \cite{huang2024high}:
\begin{align*}
  \mu[\phi^{w+1}, \phi^{w}] = \frac{\gamma}{2}(-\Delta + \beta/\varepsilon^{2})(\phi^{w+1} + \phi^{w}) + \frac{(\nu^{w+1})^m-(\nu^w)^m}{\nu^{w+1}-\nu^w}\frac{\bar E'(\phi^{w+ 1/{2}})}{m [\nu(\phi^{w+ 1/{2}})]^{m-1}},
\end{align*}
where 
$ {\phi}^{w+1/2}=(3\phi^{w} - \phi^{w-1})/2 $ is an approximation of $ (\phi^{w+1} + \phi^{w})/2 $. The relaxed DVD scheme is then defined as follows:

	\begin{equation}
	\label{system-4-sav}
	\left\{
	\begin{aligned}
	 \phi^{w+1}&= \phi^w- \Delta t   \mu[\phi^{w+1},\phi^w],\\
	\nu^{w+1}&=\nu^w+\frac{\left<\bar E'(\phi^{w+ 1/{2}}),(\phi^{w+1}-\phi^w)\right>}{m \nu^{m-1}(\phi^{w+ 1/{2}})}.
	\end{aligned}
	\right.
	\end{equation}
According to \cite{huang2024high}, we have
$$
\bar{\mathcal{E}}[\phi^{w+1},\nu^{w+1}]-\bar{\mathcal{E}}[\phi^w,\nu^w]=-\Delta t\left<\mu[\phi^{w+1},\phi^w],\mu[\phi^{w+1},\phi^w]\right>\leq 0,
$$
which corresponds to the unconditional modified energy stability of scheme \eqref{system-4-sav}. For   $m=2$, the scheme \eqref{system-4-sav} is equivalent to the SAV/CN scheme. In this work, we take $m=1$ suth that \eqref{system-4-sav} becomes a fully decoupled scheme.

Let us consider $\Omega=[0, 2\pi]^{3}$ as an example to introduce the spatial discretization. The domain is covered by a uniform mesh with the mesh size $h=2\pi/n$ and $n=2^D$. We define the grid points as $\vb{x}_{\vb{j}}:=h\vb{j}$, where $\vb{j}=(j_1 - 0.5,\,j_2 - 0.5 ,\,j_3 - 0.5)$, and ${j}_i=1,\,2,\,\cdots,\,n$. The solution of the Allen--Cahn equation at $t=t^w$ is approximated as a third-order tensor $\tensor{\phi}^w:=\left(\tensor{\phi}^w_{j_1,j_2,j_3}\right)\in\mathbb{R}^{n\times n\times n}$, where $\tensor{\phi}^w_{j_1,j_2,j_3}:=\tensor{\phi}^w_{\vb{j}}\approx \phi^w(\vb{x}_{\bm{j}})$. By discretizing the Laplace operator 
$\Delta$ using a 7-point second-order central finite difference scheme, we obtain the fully discrete system for the Allen--Cahn equation \eqref{eq:gradient-flow} as following third-order full tensor system:
\begin{align}\label{eq:matrix-vector-equation}
  \left(\tensor{I} + \tensor{L}\right)\tensor{\phi}^{w + 1} = \left(\tensor{I} -\tensor{L}\right)\tensor{\phi}^{w} -  \frac{\Delta t}{\varepsilon^{2}}\tensor{\phi}^{w + 1/2}\odot\left(\left(\tensor{\phi}^{w + 1/2}\right)^{\odot 2} -(1+\beta)\mathds{1}\right),
\end{align}
where $ \tensor{L} = \frac{\Delta t\gamma}{2}(-\Delta_{\textnormal{d}} + \beta/\varepsilon^{2}\cdot \tensor{I}) $, $\tensor{\phi}^{w+1/2}=(3\tensor{\phi}^{w}-\tensor{\phi}^{w-1})/2$, and $\mathds{1}$ is the third-order full tensor corresponding to all-ones vector ${\mathbf{1}}$, respectively. Here $\tensor{I}$ is the identity operator and  the discrete Laplace operator $\Delta_\textnormal{d}$ is given by:
$$\Delta_\textnormal{d}=\mat{D}\otimes \mat{I}\otimes\mat{I}+ \mat{I}\otimes\mat{D}\otimes\mat{I}+ \mat{I}\otimes\mat{I}\otimes \mat{D}, $$
where $\mat{I}$ is a $2^D$ order identity matrix and $\mat{D}$ is defined as:
\begin{align*}
  \mat{D} = \frac{1}{h^2}\begin{bmatrix}
    -2 & 1 & & &  1 \\
    1 & -2 &  1 & & \\
    & \ddots & \ddots & \ddots & \\
    & & 1 & -2 & 1 \\
    1 & & & 1 & -2
  \end{bmatrix}\in\mathbb{R}^{2^{D}}\times \mathbb{R}^ {  2^{D}}.
\end{align*}


Next, we employ quantized tensor train (QTT) format \cite{khoromskij2011d,KhoromskijOseledets+2011+303+322,markeeva2021qtt,kazeev2013low} to reformulate the full tensor system \eqref{eq:matrix-vector-equation} as a QTT system. Let $\mat{\Phi}^w\in\mathbb{R}^{2^{3D}}$ be a  $3D$-th order QTT tensor that approximates the full tensor $\tensor{\phi}^w$.
Since the operators $\tensor{I}\pm \tensor{L}$ are tensor products of several one-dimensional Laplacian-like operators, and inspired by \cite{kazeev2012low}, we can reshape then into: 
\begin{equation}
\begin{aligned}\label{eq:QTToperator}
  \tensor{I} \pm  \tensor{L} & = \tensor{P}^{\pm} \otimes \tensor{I} \otimes \tensor{I} + \tensor{I} \otimes \tensor{P}^{\pm} \otimes \tensor{I} + \tensor{I} \otimes \tensor{I} \otimes \tensor{P}^{\pm} \\
  & = \begin{bmatrix} \tensor{I} & \tensor{P}^{\pm} \end{bmatrix} \bowtie \begin{bmatrix}
    \tensor{I} & \tensor{P}^{\pm} \\ & \tensor{I}
  \end{bmatrix} \bowtie \begin{bmatrix}
     \tensor{P}^{\pm}\\ \tensor{I}
  \end{bmatrix},
\end{aligned}
\end{equation}
where the symbol $ \bowtie $ denotes the rank core product \cite[Definition 2.1]{kazeev2012low}.
In \eqref{eq:QTToperator}, the QTT matrix $\tensor{P}^{\pm} $ is given by:
\begin{align*}
  \tensor{P}^{\pm} = \begin{bmatrix}
    \mat{I} & \mat{J}' & \mat{J}
  \end{bmatrix} \bowtie
  \begin{bmatrix}
    \mat{I} & \mat{J}' & \mat{J} \\
    & \mat{J} & \mat{J}' \\
    & \mat{J} & \mat{J}'
  \end{bmatrix} \bowtie
  \begin{bmatrix}
    \mat{I} & \mat{J}' & \mat{J} \\
    & \mat{J} & \\
    & & \mat{J}'
  \end{bmatrix}^{\bowtie^{(D - 3)}} \bowtie \begin{bmatrix}
    \bar\alpha_{\pm} \mat{I} + \tilde\alpha_{\pm} (\mat{J} +  \mat{J}') \\
    \tilde \alpha_{\pm}\mat{J} \\
    \tilde\alpha_{\pm} \mat{J}'
  \end{bmatrix},
\end{align*}
where $ \mat{I} = \left[\begin{smallmatrix}
  1 & 0 \\
  0 & 1
\end{smallmatrix}\right] $, $ \mat{J} = \left[\begin{smallmatrix}
  0 & 1 \\
  0 & 0
\end{smallmatrix}\right] $, $\bar\alpha_{\pm} = \frac{1}{3}\left(1 \pm \Delta t \frac{\gamma}{2} \left(\frac{6}{h^{2}}+ \frac{\beta}{\varepsilon^{2}}\right)\right)$, and $ \tilde\alpha_{\pm}  = \mp\frac{\Delta t}{h^{2}} \frac{\gamma}{2}$, respectively. Following \cite[Lemma 5.2]{kazeev2012low}, the resulting QTT matrices $\tensor{I}\pm\tensor{L}$ have ranks 
\begin{align*}
  1, \underbrace{4, \dots, 4}_{D - 1}, 2, \underbrace
  {5, \dots, 5}_{D - 1}, 2, \underbrace{4, \dots, 4}_{D - 1}, 1.
\end{align*}
Using \eqref{eq:QTToperator}, we can rewrite the fully discrete system \eqref{eq:matrix-vector-equation} in QTT form as follows:
\begin{align}
\label{eq:QTTlinear}
\left(\tensor{I} + \tensor{L}\right)\mat{\Phi}^{n + 1} = \left(\tensor{I} - \tensor{L}\right)\mat{\Phi}^{n} -  \frac{\Delta t}{\varepsilon^{2}}{\mat{\Phi}}^{n + 1/2}\odot\left(({\mat{\Phi}}^{n + 1/2})^{\odot 2} - (1 + \beta)\mathds{1}\right).
\end{align}
We apply the DMRG solver (\verb|dmrg_solve2| function from TT-Toolbox) to solve the QTT system \eqref{eq:QTTlinear}. The Hadamard products on the right hand side of \eqref{eq:QTTlinear} are computed using the HaTT algorithm. For comparison, we also conduct simulations where the Hadamard products are computed directly and then recompressed using either TT-Rounding \cite{oseledets2011tensor} or one of the three randomized algorithms (RandOrth, OrthRand, and TwoSided) proposed in \cite{al2023randomized}. Two examples are considered, one with a solution that has relatively small ranks and the other with a solution that has relatively large ranks.

In the first example, we set the mobility as $\gamma=1$, the suitable stabilization parameter $\beta=1$, and the small parameter $\varepsilon=0.1$. The computational domain $\Omega=[0, 2\pi]^{3}$. The time step size is set as $\Delta t=0.01$ and the end time is $T=0.1$.  The initial condition is given by $\phi_0=0.2\sin(x_1)\sin(x_2)$ $\sin(x_3)$. The values of $\phi_0$ at mesh points $\bm{x}_{\bm{j}}$ are assembled into a third-order full tensor $\mat{\Phi}^0$. In general, one needs to apply TT-SVD \cite{oseledets2011tensor} or TT-cross \cite{oseledets2010tt} method to obtain a low rank approximation of $\mat{\Phi}^0$. However, in this example, by exploiting the specific structure of the initial values and drawing inspiration from \cite{oseledets2013constructive}, we directly express $\mat{\Phi}^0$ in a low rank QTT format, whose QTT ranks are bounded by 2. For further details on this expansion, we refer the reader to \cite{oseledets2013constructive}. 

To evaluate the accuracy of the QTT solver with the HaTT algorithm, we perform a simulation on a $2^7\times 2^7\times 2^7$ uniform mesh. We reshape \eqref{eq:matrix-vector-equation} into a matrix-vector system and solve it using the Preconditioned Conjugate Gradient (PCG) method \cite{eisenstat1981efficient}. The solution obtained from the PCG solver serves as the reference. The distributions of the solutions obtained by both the QTT solver and the PCG solver are shown in \cref{fig:contour}, demonstrating strong agreement between the two. These results validate the accuracy of the QTT solver with the HaTT algorithm. 
The evolution of the modified total free energy for the QTT solver with the HaTT algorithm, along with other QTT recompression algorithms, is presented in \cref{fig:energy-dissipation}. The results show that all solvers preserve the (modified) energy dissipation law. Compared to the PCG solver, the relative errors of the modified free energy obtained by QTT solvers range from $10^{-8} $ to $10^{-4} $. These findings indicate that the QTT solver with TT-Rounding is the most accurate, while the QTT solver with the HaTT algorithm and the three other randomized algorithms exhibit slightly lower accuracy due to the inherent randomness.

\begin{figure}[!htbp]
  \centering
  \includegraphics[width=0.95\textwidth]{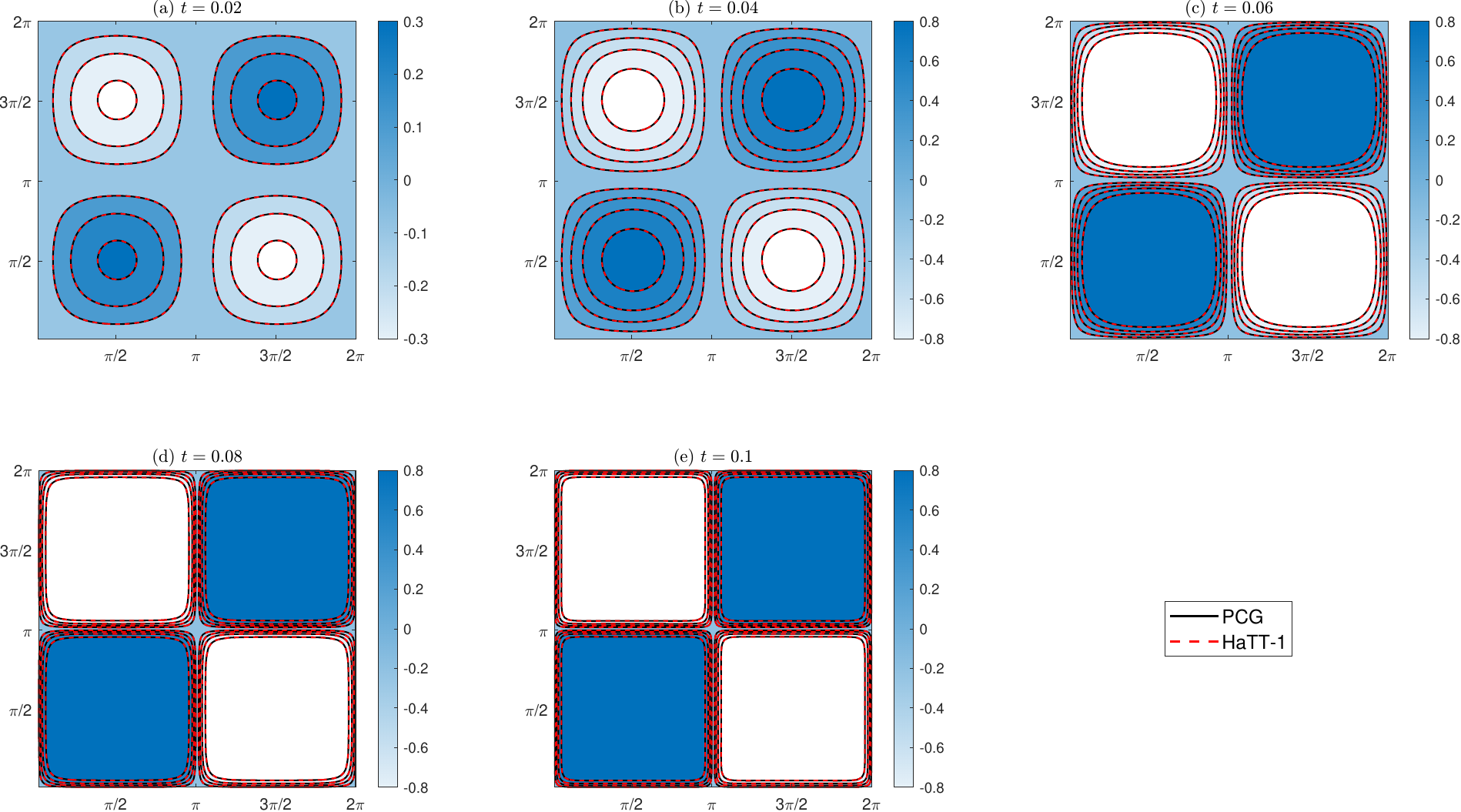}
  \caption{Contours of the solutions at different times for Example 4. }\label{fig:contour}
\end{figure}

\begin{figure}[!htpb]
  \centering
  \includegraphics[width=0.95\textwidth]{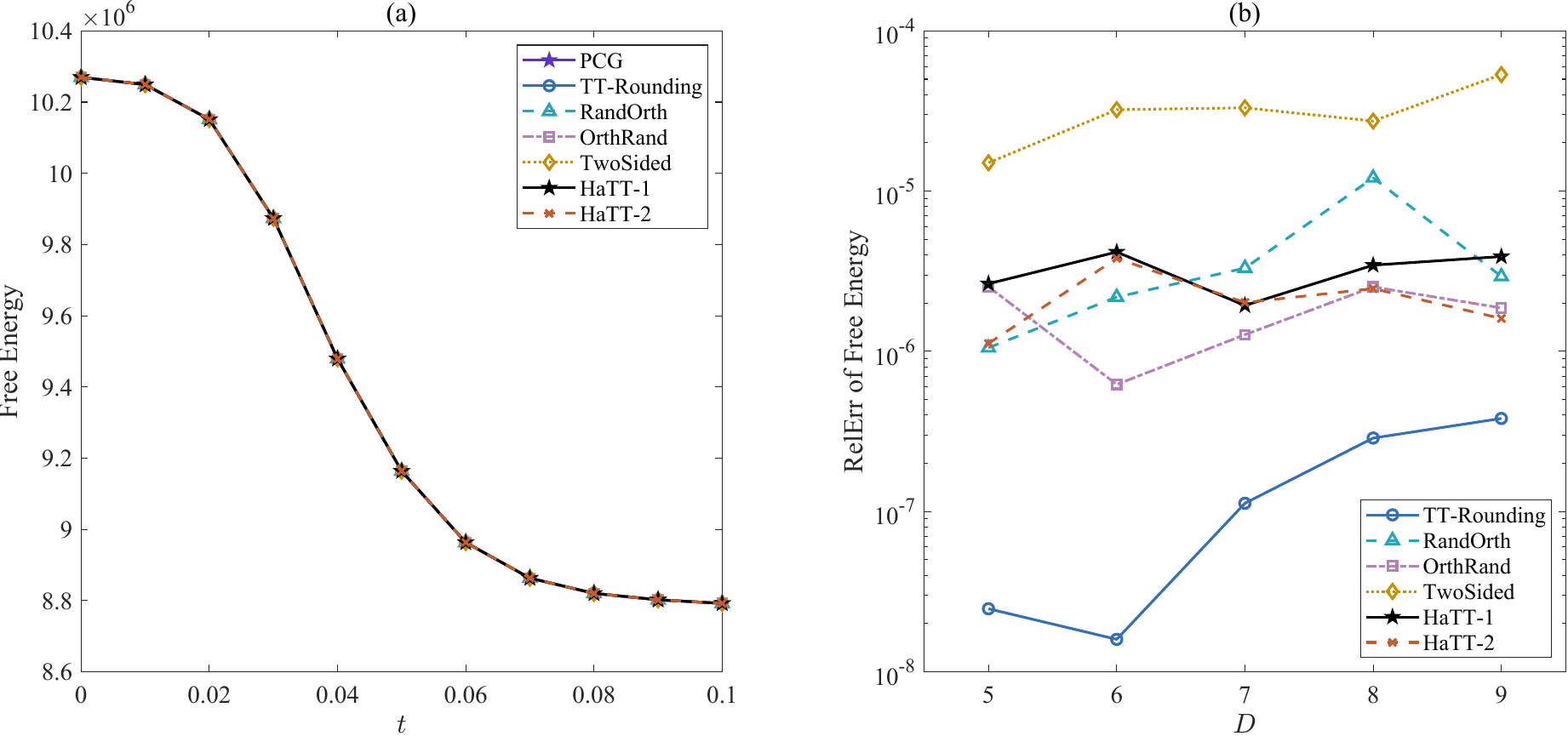}
  \caption{Example 4. (a) Modified total free energy at each time step for $ D = 7 $. (b) Relative errors of the modified free energy obtained by QTT solvers for different values of $D$.}\label{fig:energy-dissipation}
\end{figure}

To verify the complexity analysis presented in \cref{sec:complexity-analysis} and assess the efficiency of the QTT solver with the HaTT algorithm, we conduct five simulations with $D=7,\,8,\,9,\,12,$ and $15$, respectively. The computational times for these simulations with different solvers are reported in \cref{computing-time}. 
For the PCG solver, the computational complexity scales as $D2^{3D}$, meaning the computational time increases by more than a factor of 8 as $D$ increases. Simulations with $D=12$ and $D=15$ using the PCG solver are not completed due to memory limitations. In contrast, the computational cost of the QTT solver is linear in $D$, and as shown in \cref{computing-time}, the computational time increases proportionally to $D$. In all simulations, the ranks of the solution $\mat{\Phi}^n$ are largely insensitive to the mesh size $h=1/2^D$, remaining below 30. Due to the low ranks of the solutions, the QTT solver outperforms the PCG solver when $D\geq8$.

The total computational times of the QTT solvers, along with the computational times for the right hand side of \eqref{eq:QTTlinear}, are provided in \cref{computing-time}.
From these results, we observe that the performance of the QTT solver with the HaTT algorithm, the TT-Rounding algorithm, and the RandOrth algorithm are quite similar. This is because the computational time of the DMRG algorithm dominates the total computational time of the QTT solver, as shown in \cref{computing-time}. Upon examining the ranks of the solution $\mat{\Phi}^n$ and the right hand side of \eqref{eq:QTTlinear}, we find that the ranks of the right hand side are consistently twice as large as those of $\mat{\Phi}^n$. As a result, the DMRG algorithm is less efficient than the HaTT-1 algorithm or other recompression methods. 
Focusing on the computational times for the right hand side of \eqref{eq:QTTlinear}, the HaTT algorithm demonstrates a speedup ranging from $2.6\times $ to $13.7\times $ compared to the traditional TT recompression methods.

\begin{table}[!htpb]
\centering
\caption{Example 4: the total computational time of different solvers (denoted as "Total") and the computational times  for the right hand side of \eqref{eq:QTTlinear} (denoted as "Right").}\label{computing-time}
\begin{tabular}{|c|l|lllll|}
\hline
\multirow{2}{*}{Time(s)} &
  \multirow{2}{*}{Method} &
  \multicolumn{5}{c|}{$D$} \\ \cline{3-7} 
 &
   &
  \multicolumn{1}{c|}{$7$} &
  \multicolumn{1}{c|}{$8$} &
  \multicolumn{1}{c|}{$9$} &
  \multicolumn{1}{c|}{$12$} &
  \multicolumn{1}{c|}{$15$} \\ \hline
\multirow{4}{*}{Total} &
  PCG &
  \multicolumn{1}{l|}{12.04} &
  \multicolumn{1}{l|}{102.63} &
  \multicolumn{1}{l|}{1354.40} &
  \multicolumn{1}{l|}{$-$} &
  $-$
   \\ \cline{2-7} 
 &
  TTRounding &
  \multicolumn{1}{l|}{65.29} &
  \multicolumn{1}{l|}{83.29} &
  \multicolumn{1}{l|}{113.86} &
  \multicolumn{1}{l|}{353.80} &
  790.38 \\ \cline{2-7} 
 &
  RandOrth &
  \multicolumn{1}{l|}{47.46} &
  \multicolumn{1}{l|}{67.46} &
  \multicolumn{1}{l|}{82.00} &
  \multicolumn{1}{l|}{364.88} &
  880.76 \\ \cline{2-7} 
 & 
   HaTT-1 &
  \multicolumn{1}{l|}{44.23} &
  \multicolumn{1}{l|}{62.03} &
  \multicolumn{1}{l|}{77.00} &
  \multicolumn{1}{l|}{359.12} &
  874.87 \\ \hline
\multicolumn{1}{|c|}{\multirow{3}{*}{Right}} &
  TTRounding &
  \multicolumn{1}{l|}{13.98} &
  \multicolumn{1}{l|}{21.32} &
  \multicolumn{1}{l|}{29.86} &
  \multicolumn{1}{l|}{42.53} &
  54.67 \\ \cline{2-7} 
\multicolumn{1}{|l|}{} &
  RandOrth &
  \multicolumn{1}{l|}{3.32} &
  \multicolumn{1}{l|}{4.44} &
  \multicolumn{1}{l|}{4.97} &
  \multicolumn{1}{l|}{7.49} &
  9.89 \\ \cline{2-7} 
  \multicolumn{1}{|l|}{} &
  HaTT-1 &
  \multicolumn{1}{l|}{0.79} &
  \multicolumn{1}{l|}{1.01} &
  \multicolumn{1}{l|}{1.20} &
  \multicolumn{1}{l|}{1.76} &
  2.31 \\ \hline
\end{tabular}
\end{table}

Next, let us consider an example where the solutions have relatively large ranks. In this case, the computational domain is chosen as $\Omega=(-1,\,1)^2$, which is covered by a $2^D\times 2^D$ uniform mesh, with $D$ ranging from 6 to 10.  The mobility parameter is set to $\gamma=6.10351\times 10^{-5}$ and $\varepsilon=0.0078$. The time step size is set as $\Delta t=0.2$ and the end time is $T=5$. The initial condition is set as:
\begin{align*} 
  \phi_{0}(\vb{x})  = 
  \begin{cases} 
    1 & \fnorm{\vb{x}} < R_{0}, \\
    -1 & \fnorm{\vb{x}} \ge R_{0},
  \end{cases} 
\end{align*}
where $R_{0} = 100/128$. In this example, since we can not directly express $\mat{\Phi}^0$ in a low rank QTT format, we convert it into the QTT representation using TT-SVD method \cite{oseledets2011tensor}. 
We employ the relaxed DVD method with $m=1$ proposed in \cite{huang2024high}. The stabilization parameter in $\bar{E}(\phi)$ is set to $\beta = 0$. The Laplace operator $\Delta$ is discretized using a $5$-point second-order central finite difference scheme. 
The fully discrete tensor system and the corresponding QTT matrices are similar to those in the three dimensional case. Due to page limitations, we omit further details.


In this example, the maximum rank of the solution $\mat{\Phi}^n$ approaches $2^D$, nearly doubling as $D$ increases. As a result, the computational complexity of evaluating the right hand side of \eqref{eq:QTTlinear} surpasses that of the DMRG algorithm for larger values of $D$, as confirmed by the computational times reported in \cref{table:example5}. For the QTT solver with the HaTT-1 algorithm, the proportion of computational time spent on evaluating the right hand side of \eqref{eq:QTTlinear} increases from $2.3\%$ to $81.3\%$ as $D$ grows from 6 to 10. As the ranks increase, the QTT solver with TT-Rounding and the RandOrth algorithm fails for $D=9$ and $D=10$ due to memory limitations. However, leveraging its Hadamard-avoiding advantage, the QTT solver with the HaTT-1 algorithm successfully computes accurate solutions for these cases. For $D=8$, the HaTT-1 algorithm achieves speedups of  $6.34\times$ and $9.33\times$ compared to TT-Rounding and RandOrth algorithm, respectively.  

Since the maximum rank of the solution $\mat{\Phi}^n$ approaches $2^D$, the QTT solver is less efficient than the PCG solver in this example. As reported in \cite{tang2024transformed}, the solution lies on a Kolmogorov manifold with a slowly decaying width, which explains the high ranks of $\mat{\Phi}^n$ and the inefficiency of model reduction methods, such as the low-rank tensor decomposition presented in this paper. To address this issue, \cite{tang2024transformed} proposes a novel approach that transforms the equation such that the resulting solution lies on a Kolmogorov manifold with a rapidly decaying  width. We believe this approach could further enhance the performance of the QTT solver for this example, potentially enabling it to outperform the PCG solver in high-dimensional problems with fine mesh sizes. Investigating this direction is part of our future work.

\begin{table}
\centering
\caption{Example 5: the total computational time of different solvers (denoted as "Total") and the computational times  for the right hand side of \eqref{eq:QTTlinear} (denoted as "Right").}\label{table:example5}
\begin{tabular}{|c|l|lllll|}
  \hline
  \multirow{2}{*}{Time(s)} &
    \multirow{2}{*}{Method} &
    \multicolumn{5}{c|}{$D$} \\ \cline{3-7} 
   &
     &
    \multicolumn{1}{c|}{$6$} &
    \multicolumn{1}{c|}{$7$} &
    \multicolumn{1}{c|}{$8$} &
    \multicolumn{1}{c|}{$9$} &
    \multicolumn{1}{c|}{$10$} \\ \hline
  \multirow{3}{*}{Total} &  
    TTRounding &
    \multicolumn{1}{l|}{19.66} &
    \multicolumn{1}{l|}{46.35} &
    \multicolumn{1}{l|}{290.46} &
    \multicolumn{1}{l|}{$-$} &
    $-$ \\ \cline{2-7} 
   &
    RandOrth &
    \multicolumn{1}{l|}{20.86} &
    \multicolumn{1}{l|}{52.29} &
    \multicolumn{1}{l|}{404.08} &
    \multicolumn{1}{l|}{$-$} &
    $-$ \\ \cline{2-7} 
    &
    HaTT-1 &
    \multicolumn{1}{l|}{19.71} &
    \multicolumn{1}{l|}{39.68} &
    \multicolumn{1}{l|}{91.49} &
    \multicolumn{1}{l|}{288.40} &
    1666.67 \\ 
    \hline
  \multicolumn{1}{|c|}{\multirow{3}{*}{Right}} &
    TTRounding &
    \multicolumn{1}{l|}{2.07} &
    \multicolumn{1}{l|}{19.92} &
    \multicolumn{1}{l|}{237.81} &
    \multicolumn{1}{l|}{$-$} &
    $-$ \\ \cline{2-7} 
  \multicolumn{1}{|l|}{} &
    RandOrth &
    \multicolumn{1}{l|}{2.30} &
    \multicolumn{1}{l|}{24.61} &
    \multicolumn{1}{l|}{349.83} &
    \multicolumn{1}{l|}{$-$} &
    $-$ \\ \cline{2-7} 
    \multicolumn{1}{|l|}{} &
    HaTT-1 &
    \multicolumn{1}{l|}{0.45} &
    \multicolumn{1}{l|}{9.65} &
    \multicolumn{1}{l|}{37.51} &
    \multicolumn{1}{l|}{192.10} &
    1354.40 \\ \hline
  \end{tabular}
\end{table}

\section{Conclusions and limitations} \label{sec:conclusion}


We propose a Hadamard product-free TT recompression algorithm, named HaTT, for efficiently recompressing the Hadamard product of TT tensors. Utilizing the property of multiplication within the Hadamard product, the HaTT algorithm avoids the explicit representation of the Hadamard product and significantly reduces the overall computational cost of the rounding procedure compared to all existing TT recompression algorithms. Future work will focus on further analysis of the HaTT algorithm and its potential application in recompressing the Hadamard product of quantum TT tensors. {\color{black} This work primarily focuses on recompression when the target ranks are predetermined. While establishing rigorous $ \varepsilon $-stopping criteria for randomized algorithms like HaTT remains theoretically challenging, recent advances in random matrix theory suggest that probabilistic error bounds could potentially be developed for tensor decompositions. This represents an interesting direction for future theoretical work, although practical $ \varepsilon $-based stopping rules can still be implemented using empirical error monitoring.} In addition, future research will explore tolerance-controlled recompression of the Hadamard product, which is also of significant interest.

\section*{Acknowledgments}

{\color{black}We are thankful to both the reviewers for thorough reports
and useful remarks on the manuscript.}

\bibliographystyle{siamplain}
\bibliography{ref}
\end{document}